\documentclass[10pt,reqno]{amsart}
\usepackage{graphicx}
\usepackage{amsfonts}
\usepackage{amssymb}
\usepackage{amsmath}
\usepackage{amsxtra}
\usepackage{latexsym}
\usepackage{epstopdf}
\usepackage{mathrsfs}
\usepackage{esint}
\usepackage{booktabs}
\usepackage{subfigure}
\usepackage{multirow}
\usepackage{setspace}
\usepackage{enumitem}
\usepackage{comment}

\newtheorem{theorem}{Theorem}[section]
\newtheorem{lemma}[theorem]{Lemma}

\newtheorem{proposition}[theorem]{Proposition}
\newtheorem{Assumption}{Assumption}[section]

\theoremstyle{definition}

\newtheorem{algorithm}{Algorithm}[section]
\newtheorem{example}{Example}[section]

\theoremstyle{remark}
\newtheorem{remark}{Remark}[section]

\numberwithin{equation}{section}

\def\d{\delta}

\def\l{\langle}
\def\r{\rangle}

\def\la{\lambda}
\def\EE{{\mathbb E}}
\def\R{{\mathcal R}}
\def\PP{{\mathbb P}}
\def\a{\alpha}
\def\l{\langle}
\def\r{\rangle}
\def\p{\partial}
\def\ep{\varepsilon}

\title[]{\bf On convergence rates of stochastic mirror descent method under the discrepancy principle}

\author{Qinian Jin}
\address{Mathematical Sciences Institute, Australian National University, Canberra, ACT 2601, Australia}
\email{qinian.jin@anu.edu.au} \curraddr{}

\date{}

\begin{document}



\begin{abstract}
We study convergence rates of the loping stochastic mirror descent method for systems of linear and nonlinear ill-posed equations in Banach spaces with noisy data. The method combines stochastic component updates with a verification-set mechanism, yielding an \emph{a posteriori} stopping rule without repeated evaluation of the full residual. Our main result establishes order-optimal convergence rates for the actual iterate returned by this random stopping procedure under variational source conditions. The analysis exploits pathwise monotonicity of the Bregman error and stationarity after successful verification to transfer deterministic-index estimates to the random stopping index. For nonlinear problems, we introduce an exponential-weighted energy argument to control the additional perturbation terms caused by the varying derivatives.
\end{abstract}
\maketitle

\section{\bf Introduction}

We consider numerical solutions of inverse problems governed by the systems of operator equations
\begin{align}\label{sys}
F_i(x) = y_i, \qquad i = 1, \cdots, N,
\end{align}
where each $F_i: {\rm dom}(F_i)\subset X\to Y_i$ is a linear or nonlinear operator from a common Banach space $X$ to a Hilbert space $Y_i$ with domain ${\rm dom}(F_i)$. Such systems arise when measurements are collected under different experimental configurations or when a large data set is divided into smaller blocks. Throughout this paper, the exact data are assumed to be unavailable. Instead, we have noisy measurements $y_i^\d$ satisfying
\begin{align}\label{noise.2}
\|y_i^\d-y_i\| \le \d_i, \qquad i=1,\cdots,N,
\end{align}
where the noise bounds $\d_i>0$ are known. The underlying inverse problem is assumed to be ill-posed: small perturbations in the data may cause large changes in the solution, and the noisy equations may have no common solution. Regularization is therefore necessary to control noise amplification and obtain stable approximations (\cite{EHN1996}).

To incorporate prior information and select a solution, we introduce a proper, lower semi-continuous, strongly convex functional $\R: X\to (-\infty,\infty]$. Given an initial pair $x_0 \in {\rm dom}(\p \R)$ and $\xi_0 \in \p \R(x_0)$, we seek a solution $x^\dag$ satisfying
\begin{align}\label{xdag.0}
x^\dag \in \arg\min
\left\{ D_{\R}^{\xi_0}(x,x_0): x \in \bigcap_{i=1}^N {\rm dom}(F_i), F_i(x)=y_i,  i = 1, \cdots, N\right\},
\end{align}
where $\p \R$ denotes the subdifferential of $\R$ and 
\begin{align*}
D_{\R}^{\xi_0}(x,x_0) = \R(x) - \R(x_0) - \l \xi_0, x-x_0\r
\end{align*}
denotes the Bregman distance induced by \(\mathcal R\) at $x_0$ in the direction $\xi_0$. We assume that such a solution exists. Suitable strongly convex functionals can combine a stabilizing term with penalties or constraints that promote sparsity, nonnegativity, or piecewise constancy. The regularization functional thus determines both the geometry of the iteration and the preferred solution.

A natural deterministic approach is the mirror descent method applied to the data-fitting functional
\begin{align*}
\mathcal J^\d(x) = \frac{1}{2N} \sum_{i=1}^N \|F_i(x)-y_i^\d\|_{Y_i}^2.
\end{align*}
For differentiable component operators, the iteration reads
\begin{align}\label{MD}
\begin{aligned}
\xi_{n+1}^\d &= \xi_n^\d - \frac{\gamma_n}{N} \sum_{i=1}^N F_i'(x_n^\d)^* \left(F_i(x_n^\d) - y_i^\d\right),\\
x_{n+1}^\d &= \arg\min_{x\in X} \left\{\R(x) - \l\xi_{n+1}^\d, x\r\right\},
\end{aligned}
\end{align}
where $F_i'(x)$ denotes the Fr\'{e}chet derivative of $F_i$ at $x$. For a bounded linear operator $F_i$, its derivative is simply $F_i$. The method separates a dual update in $X^*$ from a primal reconstruction governed by $\R$. Strong convexity makes the latter minimization uniquely solvable and provides error control through Bregman distances. Nevertheless, noisy-data iterations must be stopped appropriately: continued iteration may eventually amplify measurement errors and deteriorate the reconstruction.

In the mirror descent method (\ref{MD}), each deterministic step requires evaluating all component residuals and the corresponding adjoint actions. This can be expensive when $N$ is large or when the evaluations involve costly forward and adjoint equations. Stochastic mirror descent method reduces the work per iteration by processing a randomly selected equation:
\begin{align*}
\begin{aligned}
\xi_{n+1}^\d &= \xi_n^\d - \gamma_n F_{i_n}'(x_n^\d)^* (F_{i_n}(x_n^\d) - y_{i_n}^\d),\\
x_{n+1}^\d &= \arg\min_{x\in X} \left\{\R(x) - \l \xi_{n+1}^\d, x\r\right\}.
\end{aligned}
\end{align*}
Here $i_n$ is sampled uniformly from $\{1,\cdots,N\}$, independently of the preceding sampling history. Before any discrepancy-based modification, the sampled direction is a conditionally unbiased estimator of the gradient of \(\mathcal J^\delta\). The method retains the regularization geometry of deterministic mirror descent while using only one component equation at each iteration.

Stochastic mirror descent method has been extensively studied in convex optimization (\cite{BT2003,NL2014,NJLS2009,ZMBBG2020}). For ill-posed inverse problems, optimization convergence must be supplemented by stability with respect to data perturbations and a suitable stopping strategy. In \cite{JL2019} the regularization property and sub-optimal convergence-rate were established for stochastic gradient descent applied to linear inverse problems under \emph{a priori} stopping rules. The corresponding analyses were developed in \cite{JZZ2020} for nonlinear ill-posed problems in Hilbert spaces. More recently, the work \cite{JLT2026} established order-optimal convergence rates for linear inverse problems in Hilbert spaces under H\"older-type source conditions with smoothness parameter $0 < \nu \le 1/2$. The results include rates in expectation and almost surely under suitable \emph{a priori} stopping rules and allow constant as well as polynomially decaying step sizes.

For linear inverse problems in Banach spaces, the stochastic mirror descent with convex regularization functionals was introduced in \cite{JLZ2023}, convergence was established under \emph{a priori} stopping rules, discrepancy-based step-size choices was investigated, and rates were derived under a benchmark source condition. The framework was extended in \cite{GFHF2025,HJLZ2025} to nonlinear ill-posed equations in Banach spaces. These developments provide a foundation for stochastic iterative regularization with penalties that reflect the structure of the desired solution.

For noisy inverse problems, however, an appropriate stopping rule is indispensable. Rate-optimal \emph{a priori} stopping rules typically depend on source parameters describing the regularity of $x^\dag$, which are generally unknown in applications. An \emph{a posteriori} rule is therefore preferable from a practical point of view, since it uses only the computed iterates, noisy measurements, and known noise bounds.

The main difficulty is to design such a stopping rule without destroying the computational advantage of stochastic updates. Evaluating the full discrepancy after every stochastic step may be prohibitively expensive. Periodic global checks reduce this burden but require a prescribed checking schedule and may allow unnecessary noisy-data updates between two checks. A natural alternative is to use the residual information already generated by the stochastic iteration.

A component discrepancy principle can be incorporated into the stochastic mirror descent method through a loping strategy, following the regularizing Landweber-Kaczmarz method of \cite{HLS2007}. The dual update becomes
\begin{align*}
\xi_{n+1}^\d = \xi_n^\d - \gamma_n a_n^\d F_{i_n}'(x_n^\d)^* (F_{i_n}(x_n^\d) - y_{i_n}^\d), \quad
a_n^\d = {\bf 1}_{\{\|F_{i_n}(x_n^\d) - y_{i_n}^\d\| > \tau \d_{i_n}\}},
\end{align*}
where $\tau>1$ and ${\bf 1}_S$ denotes the characteristic function of a set $S$. Thus, an update is omitted whenever the sampled residual is consistent with its noise level. Successful component tests must be collected at an unchanged iterate, since an active update may invalidate earlier discrepancy information.  

Finite termination and qualitative convergence for related a posteriori procedures have already been established. In \cite{HJLZ2025}, stochastic inner iterations were combined with outer steps using global residual information. In \cite{GJ2026}, a verification-set mechanism was introduced for an adaptive stochastic heavy-ball method, avoiding scheduled evaluations of the global residual. The verification set records the equations that still require verification, is reset whenever an update invalidates previously collected information, and determines termination once verification is complete. We adopt this verification principle for loping stochastic mirror descent. Our main objective is to complement these finite-termination and qualitative convergence results by establishing convergence rates under variational source conditions, a substantially more challenging task.

The random stopping index creates a substantial analytical difficulty. If $n_\d$ denotes the verification-based stopping index and
\[
\Delta_n^\delta
=D_{\mathcal R}^{\xi_n^\delta}(x^\dagger,x_n^\delta),
\]
then estimates for $\EE\left[\Delta_n^\delta\right]$ at deterministic indices cannot generally be evaluated at $n_\d$, because the stopping event $\{n_\delta=n\}$ and the reconstruction error are dependent. Loping further complicates comparisons between exact-data and noisy-data iterations: perturbing the data can change the active-update decisions and hence the subsequent iterates.

Our analysis overcomes the stopping-index difficulty through a structural property of the loping iteration. Under admissible step sizes, the Bregman error is nonincreasing along every sample path. Once verification is complete, all component discrepancy inequalities hold at the same iterate, and every subsequent update is suppressed. Consequently, the continued iteration is stationary after termination, and
\begin{align*}
\Delta_{n_\d}^\d \le \Delta_n^\d \qquad \text{for every deterministic } n\ge 0
\end{align*}
along each sample path. This pathwise comparison transfers deterministic-horizon bounds to the actual a posteriori output without decoupling the stopping index from the iterates or estimating its distribution. The deterministic comparison horizon is used only in the proof and need not be known when implementing the stopping rule.

To obtain quantitative error bounds, we impose a variational source condition. Such conditions describe solution regularity through inequalities involving the regularization functional and the data discrepancy and provide a flexible alternative to exact range-type source representations; see \cite{HKPS2007,HM2012}. Writing
\begin{align*}
f(x) = \R(x) -\l \xi_0, x\r,
\end{align*}
we consider
\begin{align*}
f(x^\dag) - f(x) \le M \|F(x)-y\|_Y^q, \qquad 0 < q \le 1
\end{align*}
on the relevant domain, where $Y = Y_1\times \cdots \times Y_N$ carries its Hilbert product norm. The exponent $q$ quantifies the regularity available for deriving rates. The case \(q=1\) includes the benchmark source condition, while \(0<q<1\) permits weaker regularity without requiring an exact range representation. The auxiliary-dual strategy developed in \cite{Jin2025} for benchmark source conditions provides an important starting point. Extending it to variational source conditions requires comparison with approximate source elements and control of their approximation defects.

The main contribution of this paper is a quantitative regularization theory establishing order-optimal convergence rates for the verification-stopped iterate produced by loping stochastic mirror descent under variational source conditions. The main ingredients and contributions can be summarized as follows.

First, we derive deterministic-horizon error estimates that are sufficiently flexible to accommodate variational source conditions and then transfer these estimates to the verification-based random stopping index. The transfer relies only on pathwise monotonicity of the Bregman error and stationarity of the algorithm after successful verification. In this way we obtain quantitative bounds for the reconstruction actually returned by the implemented \emph{a posteriori} procedure, rather than for an auxiliary deterministic iterate.

Second, for linear inverse problems we develop estimates involving an arbitrary dual comparison element $\la$ and its associated Fenchel defect. This formulation is essential when an exact source representation is unavailable. The resulting comparison estimate reduces the rate analysis to balancing the size of the comparison element against its approximation defect. Fenchel-Rockafellar duality formula then converts this balance into a primal quadratic-penalty problem, to which the variational source condition can be applied. This yields the order-optimal convergence rates
\begin{align}\label{rates}
\EE\left[D_\R^{\xi_{n_\d}^\d} (x^\dag, x_{n_\d}^\d) \right] = O(\d^q), \qquad 
\EE\left[\|x_{n_\d}^\d - x^\dag\|^2\right] = O(\d^q),
\end{align}
where 
\begin{align*}
\d := \sqrt{\d_1^2 + \cdots + \d_N^2}.
\end{align*}

Third, we extend the rate theory to nonlinear systems. The nonlinear auxiliary-dual recursion contains additional terms that couple the primal reconstruction error, the sampled residual, and the dual comparison error. These terms cannot be controlled by direct accumulation. We overcome this difficulty by introducing an exponential-weighted energy in which the squared dual comparison error is multiplied by an exponential of the Bregman error. The pathwise descent of the Bregman error then generates damping terms that absorb the nonlinear coupling. After adding a suitable multiple of the Bregman error, the remaining positive residual contribution is cancelled and a closed conditional drift estimate is obtained. Combined with a range-invariance condition, a tangential cone condition, and a local variational source condition, this leads again to the order-optimal rate (\ref{rates}). 

Although the linear and nonlinear analyses follow the same general strategy, we present them separately in order to make the structure of the argument clearer. The linear case isolates the basic dual-comparison mechanism and shows how variational source information enters the rate analysis in its simplest form. It also provides the benchmark estimates against which the nonlinear argument can be organized. The nonlinear theory then builds on this framework with the additional tools needed to accommodate the variation of the forward derivatives. This separation highlights the common core of the analysis while making transparent the genuinely new difficulties arising from nonlinearity.

The paper is organized as follows. Section 2 provides basic elements from convex analysis in Banach spaces. Section 3 treats linear inverse problems and derives deterministic-index estimates and convergence rates at the verification-based stopping index under variational source conditions. Section 4 develops the corresponding nonlinear theory, including the exponential-weighted auxiliary-dual estimates required to control nonlinear perturbations, and proves order-optimal rates for the stopped reconstruction. Finally, Section 5 provides some numerical experiments to verify theoretical convergence rates.

\section{\bf Preliminaries from convex analysis}\label{sec:convex-prelim}

Let $X$ be a real Banach space with dual space $X^*$, we use $\l \xi, x\r$ to denote the duality pairing between $\xi \in X^*$ and $x \in X$. For a proper convex function $\R: X \to (-\infty,\infty]$, we use ${\rm dom}(\R) := \{x\in X: \R(x) < \infty\}$ to denote its effective domain. The subdifferential of $\R$ at $x \in {\rm dom}(\R)$ is defined by
\begin{align*}
\p\R(x) := \left\{ \xi \in X^*: \R(z) \ge \R(x) + \l \xi, z-x\r \text{ for all } z\in X \right\}.
\end{align*}
If $\xi \in \p\R(x)$, then 
\begin{align*}
D_{\R}^{\xi}(z,x) := \R(z) - \R(x) - \l\xi, z-x\r, \qquad z \in {\rm dom}(\R)
\end{align*}
is called the Bregman distance induced by $\R$ at $x$ in the direction $\xi$. By the definition of the subdifferential, we always have
\begin{align*}
D_{\R}^{\xi}(z,x) \ge 0.
\end{align*}
In general, the Bregman distance is neither symmetric nor does it satisfy the triangle inequality.

The Legendre-Fenchel conjugate of $\R$ is the function $\R^*: X^* \to (-\infty,\infty]$ defined by
\begin{align*}
\R^*(\xi) := \sup_{x\in X} \left\{ \l\xi, x\r - \R(x)\right\}.
\end{align*}
The Fenchel-Young inequality states that
\begin{equation}\label{eq:FY}
\R(x) + \R^*(\xi) \ge \l\xi, x\r, \qquad x \in X,\,\xi \in X^*.
\end{equation}
Equality holds if and only if $\xi \in \p\R(x)$. When $\R$ is also lower semi-continuous, we have 
\begin{align*}
\xi \in \p\R(x) \Longleftrightarrow x \in \p\R^*(\xi) \Longleftrightarrow \R(x) + \R^*(\xi) = \l\xi, x\r.
\end{align*}
Given $x \in {\rm dom}(\R)$ and $\xi \in X^*$, the quantity 
\begin{align*}
\R^*(\xi) + \R(x) - \l\xi, x\r 
\end{align*}
is called a Fenchel defect, which is always nonnegative by the Fenchel-Young inequality. This defect measures the failure of $\xi$ to be a subgradient of $\R$ at $x$.

We say $\R$ is strongly convex if there exists $\sigma >0$ such that 
\begin{equation}\label{eq:SC}
\R(tx + (1-t)z) + \sigma t(1-t) \|x - z\|^2 \le t \R(x) + (1-t)\R(z)
\end{equation}
for all $x, z \in {\rm dom}(\R)$ and $0 \le t \le 1$. For strong convex function $\R$, if $\xi\in\p\R(x)$ then
\begin{align*}
\sigma \|z - x\|^2 \le D_{\R}^{\xi}(z,x), \qquad \forall z \in X. 
\end{align*}
Thus convergence in Bregman distance implies convergence in the norm of $X$.

Strong convexity of \(\mathcal R\) has an important dual consequence. If $\R: X \to (-\infty,\infty]$ is proper, lower semi-continuous and strongly convex in the sense of (\ref{eq:SC}), then $\R^*$ is Fr\'{e}chet differentiable on $X^*$ with  $\nabla \R^*$ mapping from $X^*$ to $X$ and 
\begin{align*}
\|\nabla \R^*(\xi) - \nabla \R^*(\zeta)\| \le \frac{1}{2\sigma} \|\xi - \zeta\|, \qquad \xi,\zeta \in X^*.
\end{align*}
Moreover, for any given $\xi \in X^*$ the minimization problem
\begin{align*}
\min_{x\in X} \left\{\R(x) - \l\xi, x\r\right\}
\end{align*}
has a unique solution given by $x_\xi := \nabla \R^*(\xi)$ and $\xi \in \p\R(x_\xi)$.

\section{\bf Convergence rates for linear inverse problems}\label{sect13.4} 

In this section we consider linear ill-posed inverse problems governed by the system
\begin{align}\label{linsys.B}
A_i x = y_i, \quad i = 1, \cdots, N,
\end{align}
where, for each $i \in \{1, \cdots, N\}$, the operator $A_i : X \to Y_i$ is bounded linear from a common Banach space $X$ into a Hilbert space $Y_i$. We assume that the sought solution is characterized by a proper, lower semi-continuous, and strongly convex functional $\R: X \to (-\infty, \infty]$, and that (\ref{linsys.B}) admits a solution in ${\rm dom}(\R)$. Our goal is to determine the unique solution $x^\dag$ of (\ref{linsys.B}) satisfying
\begin{align}\label{xdag.lin}
D_\R^{\xi_0}(x^\dag, x_0) = \min\left\{D_\R^{\xi_0}(x, x_0): A_i x = y_i, \, i = 1, \cdots, N\right\},
\end{align}
where $(x_0, \xi_0) \in X \times X^*$ is an initial pair with $\xi_0 \in \p \R(x_0)$. For convenience, let $Y:=Y_1\times\cdots\times Y_N$ and define $A: X \to Y$ by
\begin{align*}
Ax := (A_1 x, \cdots, A_N x), \qquad x\in X.
\end{align*}
Throughout this section, we impose the following standing assumptions.

\begin{Assumption}\label{dgm.ass1}
\begin{enumerate}[leftmargin = 0.8cm]
\item[{\rm (i)}] $X$ is a Banach space, $Y$ is a Hilbert space, and each $A_i : X \to Y_i$ is a bounded linear operator.

\item[{\rm (ii)}] $\R : X\to (-\infty, \infty]$ is proper, lower semi-continuous, and strongly convex in the sense of (\ref{eq:SC}). 
 
\item[{\rm (iii)}] The equation $Ax = y$ admits a solution in ${\rm dom}(\R)$.
\end{enumerate}
\end{Assumption}

Let $y^\d = (y_1^\d, \cdots, y_N^\d)$ be available noisy data satisfying (\ref{noise.2}). We now introduce the loping stochastic mirror descent method equipped with a verification-set stopping mechanism.

\begin{algorithm}\label{alg:SMD2}
{\it Let $\tau >1$ and $\{\gamma_n\}$ be a given sequence of positive numbers. Choose an initial pair $(x_0, \xi_0) \in X \times X^*$ with $\xi_0 \in \p \R(x_0)$. Let $x_0^\d = x_0$ and $\xi_0^\d := \xi_0$. Initialize the verification set $I_0^\d := \{1, \cdots, N\}$. For $n\ge 0$ do the following:

\begin{enumerate} [leftmargin = 0.9cm]
\item[{\rm (i)}] Randomly select an index $i_n\in \{1, \cdots, N\}$ via the uniform distribution;

\item[{\rm (ii)}] Calculate $r_n^\d := A_{i_n} x_n^\d - y_{i_n}^\d$ and set 
\begin{align*}
a_{n,i_n} := \begin{cases}
1 & \mbox{ if } \|r_n^\d\| > \tau \d_{i_n}, \\
0 & \mbox{ otherwise}; 
\end{cases}
\end{align*}

\item[{\rm (iii)}] Update $\xi_{n+1}^\d \in X^*$ and $x_{n+1}^\d \in X$ by 
\begin{align*}
\xi_{n+1}^\d = \xi_n^\d - \gamma_n a_{n, i_n} A_{i_n}^* r_n^\d, \qquad x_{n+1}^\d = \arg\min_{x\in X} \left\{\R(x) - \l \xi_{n+1}^\d, x\r\right\}
\end{align*}
and update the verification set according to
\begin{align*}
I_{n+1}^\d := \begin{cases}
I_n^\d \setminus \{i_n\}, & \mbox{ if } a_{n, i_n} = 0, \\
\{1, \cdots, N\}, & \mbox{ otherwsie}. 
\end{cases}
\end{align*}
\end{enumerate}
}
\end{algorithm}

The verification set $I_n^\d$ records the equations that still require verification. Whenever the sampled equation satisfies its discrepancy criterion, its index is removed from the verification set. By contrast, whenever an actual update is performed, previously collected verification information may no longer be valid, and the verification set is therefore reset to $\{1,\cdots,N\}$. The iteration terminates once every equation has been successfully verified.

Under Assumption \ref{dgm.ass1} one may use the analogous arguments in the proofs of \cite[Lemma 3.5]{JLZ2023} and \cite[Proposition 3.1]{HJLZ2025} to obtain the following result. 

\begin{lemma}\label{SMD.lem3}
Let $\R$ satisfy Assumption \ref{dgm.ass1}, consider Algorithm \ref{alg:SMD2} and assume that 
\begin{align}\label{SMD.41}
\tau >1, \qquad 0 < \bar \gamma := \sup_{n\ge 0} \gamma_n < \frac{4\sigma}{L} \left(1-\frac{1}{\tau}\right) 
\end{align}
holds, where $L := \max\{\|A_i\|^2: i = 1, \cdots, N\}$. Let $\Delta_n^\d:= D_\R^{\xi_n^\d} (x^\dag, x_n^\d)$. Then along any sample path there holds 
\begin{align}\label{SMD.40}
\Delta_{n+1}^\d - \Delta_n^\d \le - c_0 \gamma_n a_{n,i_n}^\d \|A_{i_n} x_n^\d - y_{i_n}^\d\|^2
\end{align}
for all $n\ge 0$, where $c_0:= 1 - \frac{1}{\tau} - \frac{\bar \gamma L}{4\sigma}>0$. If in addition 
\begin{align}\label{SMD.44}
s_n:=\sum_{k=0}^{n-1} \gamma_k \to \infty \quad \text{as }n \to \infty,
\end{align}
then the verification-based stopping index
\begin{align}\label{SMD_stopping}
n_\d := \inf\{n \ge 0: I_n^\d = \emptyset\}
\end{align}
is finite almost surely. 
\end{lemma}

According to \cite[Theorem 3.9]{GJ2026}, for the verification-based stopping index $n_\d$ determined by (\ref{SMD_stopping}), $x_{n_\d}^\d$ can be used to approximate $x^\dag$ as $\d \to 0$. Our main objective in this section is to establish quantitative convergence rates for this stopped iterate under variational source conditions on $x^\dag$.

By the formulation of Algorithm \ref{alg:SMD2}, the dual iterates satisfy $\xi_n^\d - \xi_0 \in \mbox{Ran}(A^*)$ for all $n \ge 0$, and hence there exists $\la_n^\d \in Y$ such that $\xi_n^\d - \xi_0 = A^* \la_n^\d$. To exploit the variational source conditions in the convergence rate analysis, it is necessary to make this representation explicit and to track the evolution of the auxiliary variable $\la_n^\d$ along the iterations. Careful inspection of the update rule in Algorithm \ref{alg:SMD2} leads us to introduce the sequence $\{\la_n^\d\}$ in $Y$ defined by $\la_0^\d = (0, \cdots, 0)\in Y$ and 
\begin{align}\label{SMD.42}
\la_{n+1}^\d = (\la_{n+1, 1}^\d, \cdots, \la_{n+1, N}^\d) 
\end{align}
for $n \ge 0$, where 
\begin{align*}
\la_{n+1, i}^\d = \begin{cases}
\la_{n,i_n}^\d - \gamma_n a_{n, i_n}^\d (A_{i_n} x_n^\d - y_{i_n}^\d), & \mbox{ if } i = i_n, \\[1ex]
\la_{n,i}^\d, & \mbox{ if } i \ne i_n. 
\end{cases}
\end{align*}
It is straightforward to verify that 
\begin{align*}
\xi_n^\d = \xi_0 + A^* \la_n^\d, \qquad \forall n \ge 0. 
\end{align*}

\begin{lemma}\label{SMD.lem4}
Let $\R$ satisfy Assumption \ref{dgm.ass1}. Consider the sequence $\{\la_n^\d\}$ defined by (\ref{SMD.42}). Then 
\begin{align*}
\EE\left[\|\la_{n+1}^\d - \la\|^2 + \frac{\bar \gamma}{c_0} \Delta_{n+1}^\d\right] 
& \le \EE\left[\|\la_n^\d - \la\|^2 + \frac{\bar \gamma}{c_0} \Delta_n^\d\right] + \frac{2 \gamma_n}{N} \left(\eta_\la - \EE\left[\Delta_n^\d\right]\right)  \\
& \quad \, + \frac{2 (\tau + 1) \gamma_n}{N} \d \EE\left[\| \la_n^\d - \la\|^2\right]^{1/2}
\end{align*}
for any $\la \in Y$, where $\eta_\la$ is the Fenchel defect defined by 
\begin{align}\label{Fenchel-defect}
\eta_\la := \R^*(\xi_0 + A^* \la) + \R(x^\dag) - \l \xi_0 + A^* \la, x^\dag\r.
\end{align}
\end{lemma}

\begin{proof}
Let $\la = (\la_1, \cdots, \la_N) \in Y$ be any element. By using the definition of $\la_{n+1}^\d$ and the polarization identity in Hilbert spaces we have
\begin{align}\label{SMD.21}
\|\la_{n+1}^\d - \la\|^2 &= \sum_{i\ne i_n} \|\la_{n+1,i}^\d - \la_i\|^2 + \|\la_{n+1, i_n}^\d - \la_{i_n}\|^2 \nonumber\\
& = \sum_{i\ne i_n} \|\la_{n,i}^\d - \la_i\|^2 
+ \|\la_{n, i_n}^\d - \la_{i_n} - \gamma_n a_{n, i_n}^\d (A_{i_n} x_n^\d - y_{i_n}^\d) \|^2 \nonumber \\
& = \|\la_n^\d - \la\|^2 + \gamma_n^2 a_{n, i_n}^\d \|A_{i_n} x_n^\d - y_{i_n}^\d \|^2 + J,
\end{align}
where 
\begin{align*}
J &:= - 2 \gamma_n a_{n, i_n}^\d \left\l A_{i_n} x_n^\d - y_{i_n}^\d, \la_{n,i_n}^\d - \la_{i_n} \right\r.
\end{align*}
Let 
\begin{align*}
\zeta_{n,i_n} :=  a_{n, i_n}^\d \left(A_{i_n} x_n^\d - y_{i_n}^\d\right) - \left(A_{i_n} x_n^\d - y_{i_n}\right).
\end{align*}
Then  
\begin{align*}
J & = - 2 \gamma_n \left\l \zeta_{n, i_n}, \la_{n,i_n}^\d - \la_{i_n} \right\r - 2 \gamma_n \left\l A_{i_n} x_n^\d - y_{i_n}, \la_{n,i_n}^\d - \la_{i_n} \right\r.
\end{align*}
By the definition of $a_{n,i_n}^\d$ we have 
\begin{align*}
\|\zeta_{n, i_n}\| & \le \|y_{i_n}^\d - y_{i_n}\| + (1-a_{n, i_n}^\d) \|A_{i_n} x_n^\d - y_{i_n}^\d\|
\le (\tau + 1) \d_{i_n},
\end{align*}
where we used the fact that $\|A_{i_n} x_n^\d - y_{i_n}^\d\| \le \tau \d_{i_n}$ when $a_{n,i_n}^\d = 0$. Therefore
\begin{align*}
J & \le 2 (\tau + 1) \gamma_n \d_{i_n} \| \la_{n,i_n}^\d - \la_{i_n}\| - 2 \gamma_n \left\l A_{i_n} x_n^\d - y_{i_n}, \la_{n,i_n}^\d - \la_{i_n} \right\r.
\end{align*}
Combining this with (\ref{SMD.21}) yields 
\begin{align*}
\|\la_{n+1}^\d - \la\|^2 
& \le \|\la_n^\d - \la\|^2 + \bar \gamma \gamma_n a_{n, i_n}^\d \|A_{i_n} x_n^\d - y_{i_n}^\d \|^2  \\
& \quad \, - 2 \gamma_n \left\l A_{i_n} x_n^\d - y_{i_n}, \la_{n,i_n}^\d - \la_{i_n} \right\r \\
& \quad \, + 2 (\tau + 1) \gamma_n \d_{i_n} \| \la_{n,i_n}^\d - \la_{i_n}\|.
\end{align*}
This together with (\ref{SMD.40}) in Lemma \ref{SMD.lem3} gives 
\begin{align*}
\|\la_{n+1}^\d - \la\|^2 + \frac{\bar \gamma}{c_0} \Delta_{n+1}^\d 
& \le \|\la_n^\d - \la\|^2 + \frac{\bar \gamma}{c_0} \Delta_n^\d - 2 \gamma_n \left\l A_{i_n} x_n^\d - y_{i_n}, \la_{n,i_n}^\d - \la_{i_n} \right\r \\
& \quad \, + 2 (\tau + 1) \gamma_n \d_{i_n} \| \la_{n,i_n}^\d - \la_{i_n}\|.
\end{align*}
Taking the expectation with respect to $i_n$ yields
\begin{align*}
& \EE_{i_n} \left[\|\la_{n+1}^\d - \la\|^2 + \frac{\bar \gamma}{c_0} \Delta_{n+1}^\d \right] \\
& \le \|\la_n^\d - \la\|^2 + \frac{\bar \gamma}{c_0} \Delta_n^\d -\frac{2 \gamma_n}{N} \sum_{i=1}^N \left\l A_i x_n^\d - y_i, \la_{n,i}^\d - \la_i \right\r \\
& \quad \, + \frac{2 (\tau + 1) \gamma_n}{N} \sum_{i=1}^N \d_i \| \la_{n,i}^\d - \la_i\| \\
& \le \|\la_n^\d - \la\|^2 + \frac{\bar \gamma}{c_0} \Delta_n^\d -\frac{2 \gamma_n}{N} \left\l A x_n^\d - y, \la_n^\d - \la\right\r \\
& \quad \, + \frac{2 (\tau + 1) \gamma_n}{N} \d \| \la_n^\d - \la\|.
\end{align*}
By virtue of the Fenchel-Young inequality we have $\l \xi_0 + A^* \la, x_n^\d\r \le \R^*(\xi_0 + A^* \la) + \R(x_n^\d)$. Therefore 
\begin{align}\label{smd.111}
- \left\l A x_n^\d - y, \la_n^\d - \la\right\r 
& = - \left\l A^* \la_n^\d - A^* \la, x_n^\d - x^\dag \right\r = \l\xi_n^\d - (\xi_0 + A^* \la), x^\dag - x_n^\d \r \nonumber\\
& \le \l \xi_n^\d, x^\dag - x_n^\d \r - \l \xi_0 + A^* \la, x^\dag\r + \R^*(\xi_0 + A^* \la) + \R(x_n^\d)  \nonumber \\
& = \eta_\la - \left(\R(x^\dag) - \R(x_n^\d) - \l \xi_n^\d, x^\dag - x_n^\d\r\right) \nonumber \\
& = \eta_\la - \Delta_n^\d. 
\end{align}
Consequently 
\begin{align*}
\EE_{i_n}\left[\|\la_{n+1}^\d - \la\|^2 + \frac{\bar \gamma}{c_0} \Delta_{n+1}^\d \right] 
& \le \|\la_n^\d - \la\|^2 + \frac{\bar \gamma}{c_0} \Delta_n^\d + \frac{2 \gamma_n}{N} \left(\eta_\la - \Delta_n^\d\right)  \\
& \quad \, + \frac{2 (\tau + 1) \gamma_n}{N} \d \| \la_n^\d - \la\|.
\end{align*}
Taking the full expectation gives 
\begin{align*}
\EE\left[\|\la_{n+1}^\d - \la\|^2 + \frac{\bar \gamma}{c_0} \Delta_{n+1}^\d\right] 
& \le \EE\left[\|\la_n^\d - \la\|^2 + \frac{\bar \gamma}{c_0} \Delta_n^\d\right] + \frac{2 \gamma_n}{N} \left(\eta_\la - \EE\left[\Delta_n^\d\right]\right)  \\
& \quad \, + \frac{2 (\tau + 1) \gamma_n}{N} \d \EE\left[\| \la_n^\d - \la\|\right].
\end{align*}
Using $\EE\left[\| \la_n^\d - \la\|\right] \le \EE\left[\| \la_n^\d - \la\|^2\right]^{1/2}$, we thus complete the proof.
\end{proof} 

In order to use Lemma \ref{SMD.lem4} to derive useful estimates, we need the following elementary result.

\begin{lemma}\label{Land:lem3}
Let $\{a_n\}$, $\{b_n\}$ and $\{\gamma_n\}$ be sequences of nonnegative numbers such that
$$
a_n^2 \le b_n^2 + \sum_{j=0}^n \gamma_j a_j, \quad n=0, 1,\cdots.
$$
If $\{b_n\}$ is non-decreasing, then
$$
a_n \le b_n + \sum_{j=0}^n \gamma_j, \quad n=0, 1, \cdots.
$$
\end{lemma}

\begin{proof}
This result is a slight modification of \cite[Lemma 4.6]{JLZ2023}, and it can be proved by a similar argument.
\end{proof}

Based on Lemma \ref{SMD.lem4}, we are now in a position to derive an error estimate at deterministic iteration indices. A key feature of the following result is that the estimate is formulated for an arbitrary $\la \in Y$, together with the corresponding Fenchel defect $\eta_\la$. This flexibility will be essential for incorporating variational source conditions in the subsequent analysis.

\begin{lemma}\label{SMD.lem5}
Let $\R$ satisfy Assumption \ref{dgm.ass1} hold and consider Algorithm \ref{alg:SMD2} with $\{\gamma_n\}$ satisfying (\ref{SMD.41}). Then for any $\la \in Y$ there hold
\begin{align}\label{SMD.32}
\EE\left[\Delta_n^\d\right]
\le C \left( \frac{1}{1+s_n/N} + \frac{\|\la\|^2}{1+s_n/N} + \eta_\la + \frac{s_n}{N} \d^2\right)
\end{align}
for all $n \ge 0$, where $\eta_\la$ is defined by (\ref{Fenchel-defect}), $s_n$ is defined in (\ref{SMD.44}), and $C$ is a constant depending only on $\bar \gamma$, $\sigma$, $L$, $\tau$ and $\Delta_0$. 
\end{lemma}

\begin{proof}
Recursively using the inequality established in Lemma \ref{SMD.lem4} gives 
\begin{align}\label{SMD.30}
& \EE \left[\|\la_n^\d - \la\|^2 + \frac{\bar \gamma}{c_0} \Delta_n^\d\right] 
+ \frac{2}{N} \sum_{k=0}^{n-1} \gamma_k  \EE\left[\Delta_k^\d\right] \nonumber \\
& \le \frac{\bar \gamma}{c_0} \Delta_0 + \|\la\|^2 + \frac{2}{N} \eta_\la s_n  
+ \frac{2(1+\tau)}{N} \d \sum_{k=0}^{n-1} \gamma_k \EE\left[\|\la_k^\d - \la\|^2\right]^{1/2}
\end{align}
for all $n \ge 0$. From (\ref{SMD.30}) it follows that
\begin{align*}
\EE\left[\|\la_n^\d - \la\|^2\right]
& \le \frac{\bar \gamma \Delta_0}{c_0} + \|\la\|^2 + \frac{2}{N} \eta_\la s_n  
+ \frac{2 (1+\tau)}{N} \d \sum_{k=0}^{n-1} \gamma_k \EE\left[\|\la_k^\d - \la\|^2\right]^{1/2}
\end{align*}
for all $n \ge 0$. Applying Lemma \ref{Land:lem3} and using the inequality $\sqrt{a+b}\le \sqrt{a} +\sqrt{b}$ for any $a, b\ge 0$, we obtain
\begin{align*}
\EE\left[\|\la_n^\d - \la\|^2\right]^{1/2}  
& \le \sqrt{\frac{\bar \gamma \Delta_0}{c_0} + \|\la\|^2 + \frac{2}{N} \eta_\la s_n} + \frac{2 (1+\tau)}{N} s_n \d \\
& \le \sqrt{\frac{\bar \gamma \Delta_0}{c_0} + \|\la\|^2} + \sqrt{\frac{2 \eta_\la s_n}{N}} + \frac{2(1+\tau) s_n}{N} \d.
\end{align*}
Consequently
\begin{align*}
\sum_{k=0}^{n-1} \gamma_k \EE\left[\|\la_k^\d - \la\|^2\right]^{1/2}
& \le \sqrt{\frac{\bar \gamma \Delta_0}{c_0} + \|\la\|^2} s_n
+ \sum_{k=0}^{n-1} \gamma_k \left(\sqrt{\frac{2 \eta_\la s_k}{N}} + \frac{2(1+\tau) s_k}{N} \d \right).
\end{align*}
Note that for any $\a>0$ we have 
\begin{align*}
\sum_{k=0}^{n-1} \gamma_k s_k^\a 
= \sum_{k=0}^{n-1} (s_{k+1} - s_k) s_k^\a \le \sum_{k=0}^{n-1} \int_{s_k}^{s_{k+1}} t^\a dt 
= \int_0^{s_n} t^\a dt = \frac{s_n^{1+\a}}{1+\a}.
\end{align*}
Therefore 
\begin{align*}
\sum_{k=0}^{n-1} \gamma_k \EE\left[\|\la_k^\d - \la\|^2\right]^{1/2}
& \le \sqrt{\frac{\bar \gamma \Delta_0}{c_0} + \|\la\|^2} s_n
+ \sqrt{\frac{8\eta_\la}{9N}} s_n^{3/2} + \frac{(1+\tau) s_n^2}{N} \d.
\end{align*}
Combining this with (\ref{SMD.30}) yields
\begin{align*}
\frac{\bar \gamma}{c_0} \EE\left[\Delta_n^\d\right] & +\frac{2}{N} \sum_{k=0}^{n-1} \gamma_k \EE\left[\Delta_k^\d\right] 
\le \frac{\bar \gamma \Delta_0}{c_0} + \|\la\|^2 + \frac{2\eta_\la s_n}{N} \\
& + \frac{2 (1+\tau)}{N} \d \left(\sqrt{\frac{\bar \gamma \Delta_0}{c_0} + \|\la\|^2} s_n
+ \sqrt{\frac{8\eta_\la}{9N}} s_n^{3/2} + \frac{(1+\tau) s_n^2}{N} \d\right).
\end{align*}
By using the inequalities
$$
\frac{2(1+\tau)}{N} \sqrt{\frac{\bar \gamma \Delta_0}{c_0} + \|\la\|^2} s_n \d
\le \frac{\bar \gamma \Delta_0}{c_0} + \|\la\|^2 + \frac{(1+\tau)^2 s_n^2}{N^2} \d^2 
$$
and 
\begin{align*}
\frac{2(1+\tau)}{N} \sqrt{\frac{8\eta_\la}{9N}} s_n^{3/2} \d \le \frac{8\eta_\la}{9N} s_n + \frac{(1+\tau)^2 s_n^2}{N^2} \d^2,
\end{align*}
we further obtain
\begin{align*}
\frac{\bar\gamma}{c_0} \EE\left[\Delta_n^\d\right] +\frac{2}{N} \sum_{k=0}^{n-1} \gamma_k \EE\left[\Delta_k^\d\right] 
& \le \frac{2\bar\gamma \Delta_0}{c_0} + 2\|\la\|^2 + \frac{3\eta_\la s_n}{N} + \frac{4(1+\tau)^2 s_n^2}{N^2} \d^2.
\end{align*}
According to Lemma \ref{SMD.lem3} we have $\EE\left[\Delta_k^\d\right] \ge \EE\left[\Delta_n^\d\right]$ for all $0\le k\le n-1$. Thus
\begin{align*}
\sum_{k=0}^{n-1} \gamma_k \EE\left[\Delta_k^\d\right] \ge s_n \EE\left[\Delta_n^\d\right].
\end{align*}
Therefore  
\begin{align*}
\left(\frac{\bar\gamma}{c_0} + \frac{2 s_n}{N}\right)  \EE\left[\Delta_n^\d\right] 
& \le \frac{2\bar\gamma \Delta_0}{c_0} + 2\|\la\|^2 + \frac{3\eta_\la s_n}{N} + \frac{4(1+\tau)^2 s_n^2}{N^2} \d^2.
\end{align*}
Dividing the both sides of the above equation by $\bar\gamma/c_0 + 2 s_n/N$, we obtain (\ref{SMD.32}). 
\end{proof}

We will use Lemma \ref{SMD.lem3} and Lemma \ref{SMD.lem5} to derive convergence rates for Algorithm \ref{alg:SMD2} under the variational source conditions on $x^\dag$ specified in the following assumption. 

\begin{Assumption}\label{dgm.ass2}
There exists $0 < q \le 1$ and $M > 0$ such that 
$$
f(x^\dag) - f(x) \le M \|A x - y\|^q, \quad \forall x \in \emph{dom}(\R),
$$
where $f(x) := \R(x) - \l \xi_0, x\r$. 
\end{Assumption}

\begin{remark}\label{rem:vsc}
Assumption \ref{dgm.ass2} is a H\"older-type variational source condition. Variational source conditions were introduced as a flexible alternative to the classical range-type spectral source conditions (\cite{HKPS2007}) and have subsequently been widely used in the convergence-rate analysis of linear and nonlinear inverse problems; see, for instance, \cite{F2018,HM2012,HW2017,Jin2022} and the references therein. A general variational source condition takes the form
\begin{align*}
\mathcal E^\dag(x) \le f(x) - f(x^\dag) + \varphi(\|Ax-y\|),
\end{align*}
where $\mathcal E^\dag$ is a nonnegative error functional and $\varphi$ is an index function. Assumption \ref{dgm.ass2} corresponds to the weaker case where the error term is omitted and $\varphi(t) = Mt^q$. It includes the classical range condition as a special case: if $\xi_0 + A^*w\in\partial\mathcal R(x^\dagger)$, then
\begin{align*}
f(x^\dag) - f(x) & = \R(x^\dag) - \R(x) - \l \xi_0, x^\dag - x\r \\
& \le \R(x^\dag) - \left(\R(x^\dag) + \l \xi_0 + A^* w, x - x^\dag\r\right) - \l \xi_0, x^\dag - x\r \\
& = \l w, y - A x\r \le \|w\|\,\|Ax-y\|,
\end{align*}
so Assumption~\ref{dgm.ass2} holds with $q=1$ and $M=\|w\|$. The range $0<q<1$ permits weaker source smoothness without requiring an exact range representation. This increased flexibility is one of the principal advantages of variational source conditions in the analysis of ill-posed problems in Banach spaces.
\end{remark}

We need the following version of the Fenchel-Rockafellar duality formula (see \cite[Theorem 4.4.3]{BZ2005}).

\begin{proposition}\label{dgm.prop21}
Let $X$ and $Y$ be normed spaces, let $f : X \to (-\infty, \infty]$ and $g : Y \to (-\infty, \infty]$ 
be proper, convex functions, and let $A : X \to Y$ be a bounded linear operator. If there is 
$x_0 \in \emph{dom}(f)$ such that $A x_0 \in \emph{dom}(g)$ and $g$ is continuous at
$A x_0$, then
\begin{align}\label{dgm.FRD}
\inf_{x\in X} \{f(x) + g(Ax)\} = \sup_{\eta \in Y^*} \{-f^*(A^*\eta) - g^*(-\eta) \}.
\end{align}
\end{proposition}

Now we are ready to derive the estimates on $\EE[\Delta_n^\d]$ for deterministic iteration indices under the variational source conditions. 

\begin{lemma}\label{SMD.lem6}
Let $\R$ satisfy Assumption \ref{dgm.ass1} hold and consider Algorithm \ref{alg:SMD2} with $\{\gamma_n\}$ satisfying (\ref{SMD.41}). If $x^\dag$ satisfies the variational source condition specified in Assumption \ref{dgm.ass2}, then 
\begin{align*}
\EE\left[\Delta_n^\d\right] \le C \left( \left(1+ \frac{s_n}{N}\right)^{-\frac{q}{2-q}} + \frac{s_n}{N} \d^2\right)
\end{align*}
for all integers $n \ge 0$, where $C$ is a constant depending only on $\bar \gamma$, $\sigma$, $L$, $\tau$, $q$, $M$ and $\Delta_0$. 
\end{lemma}

\begin{proof}
From Lemma \ref{SMD.lem5} it follows that  
\begin{align}\label{SMD.22}
\EE\left[\Delta_n^\d\right]
\le C \left( \frac{1}{1+s_n/N} + \inf_{\la\in Y} \left\{\frac{\|\la\|^2}{1+s_n/N} + \eta_\la\right\} + \frac{s_n}{N} \d^2\right).
\end{align}
Recalling the definition of $\eta_\la$ and using the Fenchel-Rockafellar duality formula, we have 
\begin{align*}
& \inf_{\la\in Y} \left\{\frac{\|\la\|^2}{1+s_n/N} + \eta_\la\right\}  \\
& = \R(x^\dag) - \l \xi_0, x^\dag\r - \sup_{\la \in Y} \left\{-\R^*(\xi_0 + A^* \la) + \l \la, y\r - \frac{\|\la\|^2}{1+s_n/N}\right\} \\
& = \R(x^\dag) - \l \xi_0, x^\dag\r - \inf_{x \in X} \left\{\R(x) - \l \xi_0, x\r + \frac{1+s_n/N}{4} \|A x - y\|^2\right\} \\
& = \sup_{x \in X} \left\{f(x^\dag) - f(x) - \frac{1+s_n/N}{4} \|A x - y\|^2\right\},
\end{align*}
where $f(x) := \R(x) - \l\xi_0, x\r$. By using the variational source condition specified in Assumption \ref{dgm.ass2}, we then obtain 
\begin{align*}
\inf_{\la\in Y} \left\{\frac{\|\la\|^2}{1+s_n/N} + \eta_\la\right\}
& \le \sup_{x \in X} \left\{M \|A x - y\|^q - \frac{1+s_n/N}{4} \|A x - y\|^2\right\} \\
& \le \sup_{t \ge 0} \left\{M t^q - \frac{1+s_n/N}{4} t^2\right\} \\
& = \left(1 -\frac{q}{2}\right) M (2qM)^{\frac{q}{2-q}} \left(1 + \frac{s_n}{N}\right)^{-\frac{q}{2-q}}. 
\end{align*}
Combining this with (\ref{SMD.22}) and using $q/(2-q) \le 1$ completes the proof. 
\end{proof}

Finally we are in a position to show the convergence rate result concerning Algorithm \ref{alg:SMD2} for the stopping index $n_\d$ output by the {\it a posteriori} rule (\ref{SMD_stopping}).

\begin{theorem}\label{SMD.thm5}
Let $\R$ satisfy Assumption \ref{dgm.ass1}, and consider Algorithm \ref{alg:SMD2} with $\gamma_n$ satisfying (\ref{SMD.41}) and (\ref{SMD.44}). Let $n_\d$ be the stopping index defined by (\ref{SMD_stopping}), which is finite almost surely. If $x^\dag$ satisfies the variational source condition specified in Assumption \ref{dgm.ass2}, then 
\begin{align*}
\EE\left[D_\R^{\xi_{n_\d}^\d}(x^\dag,x_{n_\d}^\d)\right] \le C\d^q
\quad \mbox{ and } \quad \EE\left[\|x_{n_\d}^\d-x^\dag\|^2\right] \le C\d^q,
\end{align*}
where $C$ is a constant depending only on $\bar \gamma$, $\sigma$, $L$, $\tau$, $q$, $M$ and $\Delta_0$. 
\end{theorem}

\begin{proof}
Let $\Delta_n^\d:=D_\R^{\xi_n^\d}(x^\dag,x_n^\d)$. We first claim that, for every deterministic integer $n \ge 0$,
\begin{align}\label{SMD.43}
\Delta_{n_\d}^\d \le \Delta_n^\d
\end{align}
along every sample path in the event $\{n_\d < \infty\}$. Indeed, fix a sample path. If $n\le n_\d$, then (\ref{SMD.43}) follows directly from the monotonicity established in Lemma \ref{SMD.lem3}. If $n>n_\d$, then, by the definition of $n_\d$, we have
\begin{align*}
\|A_i x_{n_\d}^\d - y_i^\d\| \le \tau\d_i, \qquad i = 1, \cdots, N.
\end{align*}
Thus, by the formulation of Algorithm \ref{alg:SMD2}, no further updates are performed after $n_\d$, and hence $\xi_n^\d=\xi_{n_\d}^\d$ and $x_n^\d=x_{n_\d}^\d$. Consequently $\Delta_n^\d=\Delta_{n_\d}^\d$ and (\ref{SMD.43}) follows.

Now, since $\PP(n_\d < \infty) =1$, by taking the full expectation on the equation (\ref{SMD.43}), we obtain for every deterministic integer $n \ge 0$ that 
\begin{align*}
\EE\left[\Delta_{n_\d}^\d \right] \le \EE\left[\Delta_n^\d\right].
\end{align*}
Therefore, it follows from Lemma \ref{SMD.lem6} that  
\begin{align*}
\EE\left[\Delta_{n_\d}^\d \right] \le C\left(\left(1 + \frac{s_n}{N}\right)^{-\frac{q}{2-q}} + \frac{1}{N} s_n \d^2\right)
\end{align*}
for all deterministic integer $n\ge 0$. Choose the deterministic integer $\hat n_\d$ such that 
$$
\frac{1}{N} s_{\hat n_\d} \le \d^{q-2} < \frac{1}{N} s_{\hat n_\d + 1}. 
$$
Because of (\ref{SMD.44}), this $\hat n_\d$ is well-defined. Moreover, $\d^{q-2} \le \frac{1}{N} (s_{\hat n_\d} +\bar \gamma) \le C(1 + s_{\hat n_\d}/N)$. Thus  
\begin{align*}
\EE\left[\Delta_{n_\d}^\d \right]  
\le C\left( \left(1+ \frac{s_{\hat n_\d}}{N}\right)^{-\frac{q}{2-q}} + \frac{1}{N} s_{\hat n_\d} \d^2\right) \le C \d^q.
\end{align*}
This together with the strong convexity of $\R$ gives $\EE\left[\|x_{n_\d}^\d-x^\dag\|^2\right] \le C\d^q$ and the proof is complete. 
\end{proof}

\begin{remark}
When $X$ is a Hilbert space and there exist $0< \nu \le 1/2$ and $\omega \in X$ such that 
\begin{align}\label{ssc}
\xi_0 + (A^*A)^\nu \omega \in \p \R(x^\dag),
\end{align}
by the interpolation inequality, it is easy to see that the variational source condition in Assumption \ref{dgm.ass2} holds with $q = \frac{4\nu}{1+\nu}$ and $M = \frac{1+\nu}{2} \left(\frac{1-\nu}{\sigma}\right)^{\frac{1-\nu}{1+\nu}} \|\omega\|^{\frac{2}{1+\nu}}$; see \cite[Remark 3.3]{JW2023}. Therefore, under (\ref{ssc}), for Algorithm \ref{alg:SMD2} and the stopping index $n_\d$ determined by (\ref{SMD_stopping}), there holds 
$$
\EE\left[\|x_{n_\d}^\d - x^\dag\|^2\right] = O(\d^{\frac{4\nu}{1+2\nu}}).
$$
\end{remark}

\begin{remark}
A careful inspection of the proof shows that the strategy developed in this section yields an alternative derivation of the convergence rates for the stochastic gradient method established in \cite{JLT2026}. Moreover, it complements those results in \cite{JLT2026} by providing an order-optimal \emph{a posteriori} stopping rule.
\end{remark}

\section{\bf Convergence rates for nonlinear inverse problems}\label{sect13.5}

In this section, we derive convergence rates for the stochastic mirror descent method applied to the nonlinear ill-posed system (\ref{sys}) with noisy data $y^\d := (y_1^\d, \cdots, y_N^\d)$ satisfying (\ref{noise.2}). For each $i = 1, \cdots, N$, the operator $F_i: X\to Y_i$ is assumed to be Fr\'{e}chet differentiable from a common Banach space $X$ into a Hilbert space $Y_i$. We write  
$$
F(x) := (F_1(x), \cdots, F_N(x)), \quad x \in X. 
$$
Our aim is to recover a solution whose structural properties are encoded by a proper, lower semi-continuous, and strongly convex functional $\R: X \to (-\infty, \infty]$. Given an initial pair $(x_0, \xi_0) \in X \times X^*$ with $\xi_0 \in \p \R(x_0)$, we seek a solution $x^\dag$ of (\ref{sys}) satisfying (\ref{xdag.0}).

We impose the following standard conditions on the regularization functional $\R$ and the forward operators $F_i$, $i =1, \cdots, N$. 

\begin{Assumption}\label{ass2}
\begin{enumerate}[leftmargin = 0.9cm]
\item[\emph{(i)}] $X$ is a Banach space, $Y_i$ is a Hilbert space for each $i$.

\item[\emph{(ii)}] $\R: X \to (-\infty, \infty]$ is proper, lower semi-continuous, and strongly convex in the sense of (\ref{eq:SC}). 

\item[\emph{(iii)}] There exists $\rho>0$ such that $B_{2\rho}(x_0)\subset \emph{dom}(F)$ and (\ref{xdag.0}) has a solution $\tilde x$ such that $D_\R^{\xi_0} (\tilde x, x_0) < \sigma \rho^2$.

\item[\emph{(iv)}] Each operator $F_i$ is Fr\'{e}chet differentiable over $B_{2\rho}(x_0)$, and there exists $0\le \eta < 1$ such that 
\begin{align*}
\|F_i(x) - F_i(x^\dag) - F_i'(x^\dag) (x - x^\dag)\| \le \eta \|F_i(x) - F_i(x^\dag)\|
\end{align*}
for all $x \in B_{2\rho}(x_0)$ and $i = 1, \cdots, N$.

\item[\emph{(v)}] There is $\kappa \ge 0$ such that for each $i\in \{1, \cdots, N\}$ and each $x \in B_{2\rho}(x_0)$, there exists a bounded linear operator $Q_{i,x}: Y_i \to Y_i$ such that 
\begin{align*}
F_i'(x) = Q_{i,x} F_i'(x^\dag) \quad \mbox{and} \quad \|I - Q_{i,x}\| \le \kappa \|x-x^\dag\|.
\end{align*}
\end{enumerate}
\end{Assumption}

Condition (v) is a range-invariance condition, first introduced in \cite{HNS95}. It implies the tangential cone condition in (iv) provided $\kappa\rho$ is sufficiently small; see, for instance, \cite{HNS95,Jin2025}. We nevertheless retain (iv) explicitly for convenience in the subsequent analysis. Let
\begin{align}\label{SMD.L}
L := \max_{i= 1, \cdots, N} \sup_{x\in B_{2\rho}(x_0)} \|F_i'(x)\|^2.
\end{align}
From condition (v) of Assumption \ref{ass2} it is easy to see that $L < \infty$.  

We consider the following loping stochastic mirror descent method equipped with a verification-set mechanism.

\begin{algorithm}\label{alg:SMD0}
{\it Let $\tau >1$ and$\{\gamma_n\}$ be a given sequence of positive numbers. Pick an initial pair $(x_0, \xi_0) \in X \times X^*$ with $x_0 \in {\rm dom}(F)$ and $\xi_0 \in \p \R(x_0)$. Let $x_0^\d = x_0$ and $\xi_0^\d = \xi_0$. Set $I_0^\d := \{1, \cdots, N\}$. For $n \ge 0$ do the following:
\begin{enumerate}[leftmargin = 0.9cm]
\item[{\rm (i)}] Pick $i_n \in \{1, \cdots, N\}$ randomly via the uniform distribution;

\item[{\rm (ii)}] Set 
\begin{align*}
a_{n,i_n}^\d = \begin{cases}
1, & \mbox{ if } \|F_{i_n}(x_n^\d) - y_{i_n}^\d\| > \tau \d_{i_n}, \\
0, & \mbox{ otherwise}; 
\end{cases}
\end{align*}

\item[{\rm (iii)}] Update $\xi_{n+1}^\d$ and $x_{n+1}^\d$ by 
\begin{align*}
\xi_{n+1}^\d & = \xi_n^\d - \gamma_n a_{n,i_n}^\d  F_{i_n}'(x_n^\d)^* (F_{i_n}(x_n^\d) - y_{i_n}^\d), \\
x_{n+1}^\d & = \arg\min_{x\in X} \left\{\R(x) - \l \xi_{n+1}^\d, x\r\right\}
\end{align*}
and set 
\begin{align*}
I_{n+1}^\d = \begin{cases}
I_n^\d \setminus \{i_n\}, & \mbox{ if } a_{n,i_n} = 0,\\
\{1, \cdots, N\}, & \mbox{ otherwise}.
\end{cases}
\end{align*}
\end{enumerate}
}
\end{algorithm}

Based on Assumption \ref{ass2} (i)-(iv), we may follow the similar argument in \cite{HJLZ2025} to establish the following result, which demonstrates that Algorithm \ref{alg:SMD0} is well-defined and the Bregman distance $D_\R^{\xi_n^\d}(x^\dag, x_n^\d)$ is monotonically decreasing. 

\begin{lemma}\label{Land:lem6}
Let the conditions (i)-(iv) of Assumption \ref{ass2} hold with $L<\infty$ for $L$ defined by (\ref{SMD.L}). Consider Algorithm \ref{alg:SMD0}. Let $\bar \gamma := \sup_n \gamma_n$ and assume that
\begin{align}\label{sys.2}
\tau > \frac{1+\eta}{1-\eta} \quad \mbox{ and } \quad c_1 := 1- \eta - \frac{1+\eta}{\tau} - \frac{\bar \gamma L}{4 \sigma} >0. 
\end{align}
Then $x_n^\d \in B_{2\rho}(x_0)$ and 
\begin{align*}
\Delta_{n+1}^\d - \Delta_n^\d \le - c_1 \gamma_n a_{n, i_n}^\d \|F_{i_n}(x_n^\d) - y_{i_n}^\d\|^2
\end{align*}
for all integers $n \ge 0$, where $\Delta_n^\d:= D_\R^{\xi_n^\d}(x^\dag, x_n^\d)$. If, in addition, $\{\gamma_n\}$ satisfies (\ref{SMD.44}), then the stopping index $n_\d$ defined by 
\begin{align}\label{LSMD-stop}
n_\d := \inf\{n \ge 0: I_n^\d = \emptyset\}
\end{align}
is finite almost surely. 
\end{lemma}

To derive convergence rates of Algorithm \ref{alg:SMD0} under variational source conditions to be specified later, we note that, under Assumption \ref{ass2} (v), the updating formula for $\xi_{n+1}^\d$ can be written as 
$$
\xi_{n+1}^\d = \xi_n^\d - \gamma_n a_{n,i_n}^\d A_{i_n}^* Q_{i_n,x_n^\d}^* (F_{i_n}(x_n^\d) - y_{i_n}^\d),
$$
where we set $A_i := F_i'(x^\dag)$ for $i = 1, \cdots, N$. This suggests introducing a sequence $\{\la_n^\d\}$ in $Y$, defined as follows. Let $\la_0^\d = (0, \cdots, 0) \in Y$ and, for $n \ge 0$, let $\la_{n+1}^\d = (\la_{n+1,1}^\d, \cdots, \la_{n+1,N}^\d) \in Y$ with 
\begin{align*} 
\la_{n+1, i}^\d = \begin{cases} 
\la_{n, i}^\d, & \mbox{ if } i \ne i_n, \\ 
\la_{n, i_n}^\d - \gamma_n a_{n,i_n}^\d Q_{i_n,x_n^\d}^*(F_{i_n}(x_n^\d) - y_{i_n}^\d), & \mbox{ if } i = i_n.
\end{cases}
\end{align*}
By induction on $n$, we have $\xi_n^\d = \xi_0 + A^* \la_n^\d$ for all $n \ge 0$, where $A: X \to Y_1\times \cdots \times Y_N$ is defined by $A x := (A_1 x, \cdots, A_N x)$ for all $x \in X$.

We can establish the following estimate for the sequence $\{\la_n^\d\}$, which is analogous to \cite[Lemma 4.2]{Jin2025}. The main difference is that the present result is derived directly for noisy data and incorporates an arbitrary element $\la \in Y$, a feature that will be important in the subsequent rate analysis.

\begin{lemma}\label{recursion}
Let Assumption \ref{ass2} hold. Consider Algorithm \ref{alg:SMD0} and set 
\begin{align*}
r_{n,i}^\d := \|F_i(x_n^\d) - y_i^\d\|, \quad i = 1, \cdots, N.
\end{align*}
Then, for any $\la = (\la_1, \cdots, \la_N) \in Y$,  there holds 
\begin{align*}
\|\la_{n+1}^\d - \la\|^2 
\le{}& \|\la_n^\d - \la\|^2 + (1+\kappa \rho)^2 \gamma_n^2 a_{n, i_n}^\d \|r_{n, i_n}^\d\|^2 \\
&  - 2\gamma_n a_{n,i_n}^\d \left\l A_{i_n} (x_n^\d - x^\dag), \la_{n,i_n}^\d - \la_{i_n}\right\r \\
& + (3+\eta) \kappa \gamma_n a_{n,i_n}^\d \|x_n^\d - x^\dag\| \|r_{n,i_n}^\d\| \|\la_{n,i_n}^\d - \la_{i_n}\| \\
& + (2+(1+\eta) \kappa \rho) \gamma_n a_{n,i_n}^\d \d_{i_n} \|\la_{n,i_n}^\d - \la_{i_n}\|.
\end{align*}
\end{lemma}

\begin{proof}
By the definition of $\la_{n+1}^\d$ and the polarization identity in Hilbert spaces we have 
\begin{align*}
\|\la_{n+1}^\d - \la\|^2 
& = \|\la_n^\d - \la\|^2 + \gamma_n^2 a_{n,i_n}^\d \|Q_{i_n, x_n^\d}^*r_{n,i_n}^\d\|^2  \\
& \quad \, -2 \gamma_n a_{n,i_n}^\d \l r_{n,i_n}^\d, Q_{i_n, x_n^\d} (\la_{n,i_n}^\d - \la_{i_n})\r.
\end{align*}
Note that 
\begin{align*}
r_{n,i_n}^\d = A_{i_n} (x_n^\d - x^\dag) + \left(F_{i_n}(x_n^\d)-y_{i_n}- A_{i_n}(x_n^\d - x^\dag)\right) + \left(y_{i_n} - y_{i_n}^\d\right). 
\end{align*}
Therefore   
\begin{align*}
\left\l r_{n,i_n}^\d, Q_{i_n, x_n^\d} (\la_{n,i_n}^\d - \la_{i_n})\right\r
={}& \l A_{i_n}(x_n^\d - x^\dag), \la_{n,i_n}^\d - \la_{i_n} \r \\
& + \l F_{i_n}(x_n^\d)-y_{i_n}- A_{i_n}(x_n^\d - x^\dag), \la_{n,i_n}^\d - \la_{i_n} \r \\
& + \l y_{i_n} - y_{i_n}^\d, \la_{n,i_n}^\d - \la_{i_n} \r \\
& + \l r_{n,i_n}^\d, (Q_{i_n, x_n^\d} - I) (\la_{n,i_n}^\d - \la_{i_n} \r.
\end{align*}
Consequently
\begin{align}\label{loping_smd.2}
\|\la_{n+1}^\d - \la\|^2 
\le{}& \|\la_n^\d - \la\|^2 + \gamma_n^2 a_{n,i_n}^\d \|Q_{i_n, x_n^\d}\|^2\|r_{n,i_n}^\d\|^2 \nonumber \\
& -2\gamma_n a_{n,i_n}^\d \l A_{i_n}(x_n^\d - x^\dag), \la_{n,i_n}^\d - \la_{i_n} \r \nonumber \\
& + 2 \gamma_n a_{n,i_n}^\d \|F_{i_n}(x_n^\d)-y_{i_n}-A_{i_n}(x_n^\d- x^\dag)\| \|\la_{n, i_n}^\d - \la_{i_n}\| \nonumber \\
&+2\gamma_n a_{n,i_n}^\d \left(\d_{i_n} + \|Q_{i_n, x_n^\d} - I\| \|r_{n,i_n}^\d\|\right) \|\la_{n,i_n}^\d - \la_{i_n}\|.
\end{align}
Applying the range-invariance condition (v) and then the tangential cone condition (iv) in Assumption \ref{ass2} gives
\begin{align*}
\|F_{i_n}(x_n^\d)-y_{i_n}-A_{i_n}(x_n^\d- x^\dag)\| 
& \le \frac{1}{2} \kappa \|x_n^\d - x^\dag\| \|A_{i_n} (x_n^\d-x^\dag)\| \\
& \le \frac{1}{2} \kappa (1+\eta) \|x_n^\d-x^\dag\| (\|r_{n,i_n}^\d\| + \d_{i_n}).
\end{align*}
From Lemma \ref{Land:lem6} and Assumption \ref{ass2} (ii)-(iii) it follows that 
\begin{align*}
\sigma \|x_n^\d - x^\dag\|^2 \le \Delta_n^\d \le D_\R^{\xi_0}(x^\dag, x_0) < \sigma \rho^2
\end{align*}
which implies $\|x_n^\d - x^\dag\| < \rho$. Thus 
\begin{align*}
\|Q_{i_n, x_n^\d}-I\|\le \kappa \|x_n^\d - x^\dag\|,\quad \|Q_{i_n, x_n^\d}\|\le 1 + \kappa \rho.
\end{align*}
Combining these estimates with (\ref{loping_smd.2}) shows the desired result.
\end{proof}

Because of the nonlinearity of the forward operators, the right-hand side of the inequality established in Lemma \ref{recursion} contains the term
\begin{align*} 
(3+\eta)\kappa \gamma_n a_{n,i_n}^\d \|x_n^\d - x^\dag\| \|r_{n,i_n}^\d\| \|\la_{n,i_n}^\d - \la_{i_n}\|,
\end{align*}
which couples the primal error, the sampled residual, and the dual comparison error. A direct accumulation of such terms does not yield the bounds required for the rate analysis. This difficulty is one of the main obstacles that prevented the analysis in \cite{Jin2025} from being extended to quantitative error estimates in the noisy-data setting.

We overcome this difficulty by introducing a weighted energy argument. More precisely, we weight the squared dual error by an exponential factor involving the Bregman error. The descent property of the Bregman error then generates additional damping terms, which can be used to absorb the nonlinear coupling above. This mechanism is made precise in the following result.

\begin{lemma}\label{loping_smd.lem3}
Let Assumption \ref{ass2} hold. Consider Algorithm \ref{alg:SMD0} with $\{\gamma_n\}$ satisfying (\ref{sys.2}). Let 
$$
\Delta_n^\d := D_\R^{\xi_n^\d} (x^\dag, x_n^\d), \qquad \Phi_n := e^{\vartheta \Delta_n^\d}
$$ 
with $\vartheta := (3+\eta)^2 \kappa^2/(2c_0\sigma)$. For every $\la \in Y$ let $\eta_\la$ be the Fenchel defect defined in (\ref{Fenchel-defect}). Then there holds
\begin{align}\label{weighteddrift}
\EE_{i_n} \left[\Phi_{n+1} \|\la_{n+1}^\d - \la\|^2 \right] 
& \le \Phi_n \|\la_n^\d - \la\|^2 - \frac{\gamma_n}{N} \Phi_n \Delta_n^\d + \frac{2\gamma_n}{N} \Phi_n \eta_\la \nonumber \\
&  + \frac{C_2 \gamma_n^2}{N} \Phi_n \sum_{i=1}^N a_{n,i}^\d \|r_{n,i}^\d\|^2 + \frac{C_3 \gamma_n}{N} \d \Phi_n \|\la_n^\d - \la\|, 
\end{align}
where 
\begin{align*} 
C_2 := \max\left\{(1+\kappa\rho)^2, 2(1+\eta)^2 (1+1/\tau)^2\right\}, \quad C_3 := 2 + (1+\eta) (2 + 2\tau + \kappa \rho).
\end{align*}
\end{lemma}

\begin{proof}
According to Lemma \ref{Land:lem6} we have 
\begin{align*}
\Delta_{n+1}^\d \le \Delta_n^\d - c_0 \gamma_n a_{n,i_n}^\d \|r_{n,i_n}^\d\|^2.
\end{align*}
By introducing 
\begin{align*}
w_{n,i}^\d := e^{-\vartheta c_0\gamma_n a_{n,i}^\d \|r_{n,i}^\d\|^2}, \quad i = 1, \cdots, N, 
\end{align*}
we have 
\begin{align*}
\Phi_{n+1} = e^{\vartheta \Delta_{n+1}^\d} \le \Phi_n e^{-\vartheta c_0\gamma_n a_{n,i_n}^\d \|r_{n,i_n}^\d\|^2}
=\Phi_n w_{n,i_n}^\d.
\end{align*}
Therefore, by virtue of Lemma \ref{recursion}, we obtain 
\begin{align*}
\Phi_{n+1} \|\la_{n+1}^\d - \la\|^2 
\le{}& \Phi_n w_{n,i_n}^\d \|\la_{n+1}^\d - \la\|^2 \\
\le{}& \Phi_n \|\la_n^\d - \la\|^2 w_{n,i_n}^\d + (1+\kappa\rho)^2 \gamma_n^2 \Phi_n w_{n,i_n}^\d a_{n,i_n}^\d \|r_{n,i_n}^\d\|^2 \\
& - 2\gamma_n \Phi_n w_{n,i_n}^\d a_{n,i_n}^\d \left\l A_{i_n} (x_n^\d-x^\dag), \la_{n,i_n}^\d - \la_{i_n}\right\r \\
& + (3+\eta) \kappa \gamma_n \Phi_n w_{n,i_n}^\d a_{n,i_n}^\d \|x_n^\d - x^\dag\| \|r_{n,i_n}^\d\| \|\la_{n,i_n}^\d - \la_{i_n}\| \\
& + \left(2+(1+\eta)\kappa \rho\right) \gamma_n \Phi_n w_{n,i_n}^\d a_{n,i_n}^\d \d_{i_n} \|\la_{n,i_n}^\d - \la_{i_n}\|.
\end{align*}
Taking the expectation with respect to $i_n$ gives 
\begin{align}\label{loping_smd.5}
& \EE_{i_n} \left[\Phi_{n+1} \|\la_{n+1}^\d - \la\|^2\right] \nonumber \\
& \le \frac{\Phi_n \|\la_n^\d - \la\|}{N} \sum_{i=1}^N w_{n,i}^\d + \frac{(1+\kappa\rho)^2 \gamma_n^2 \Phi_n}{N} \sum_{i=1}^N w_{n,i}^\d a_{n,i}^\d \|r_{n,i}^\d\|^2 \nonumber \\
& \quad \, + \frac{(2+(1+\eta)\kappa \rho) \gamma_n \Phi_n}{N} \sum_{i=1}^N w_{n,i}^\d a_{n,i}^\d \d_i \|\la_{n,i}^\d - \la_i\| + I_1 + I_2,
\end{align}
where 
\begin{align*}
I_1 &:= - \frac{2\gamma_n \Phi_n}{N} \sum_{i=1}^N w_{n,i}^\d a_{n,i}^\d \left\l A_i (x_n^\d-x^\dag), \la_{n,i}^\d - \la_i \right\r, \\
I_2 &:= \frac{(3+\eta)\kappa \gamma_n \Phi_n}{N} \|x_n^\d - x^\dag\| \sum_{i=1}^N w_{n,i}^\d a_{n,i}^\d  \|r_{n,i}^\d\| \|\la_{n,i}^\d - \la_i\|.
\end{align*}
By writing 
$$
w_{n,i}^\d a_{n,i}^\d = 1 - (1 - a_{n,i}^\d) - a_{n,i}^\d (1 - w_{n,i}^\d),
$$
we obtain 
\begin{align*}
I_1 ={}& - \frac{2\gamma_n \Phi_n}{N} \left\l A (x_n^\d - x^\dag), \la_n^\d - \la \right\r  \\
& - \frac{2\gamma_n \Phi_n}{N} \sum_{i=1}^N (1 - a_{n,i}^\d) \left\l A_i (x_n^\d - x^\dag), \la_{n,i}^\d - \la_i \right\r \\
& - \frac{2\gamma_n \Phi_n}{N} \sum_{i=1}^N a_{n,i}^\d (1 - w_{n,i}^\d) \left\l A_i (x_n^\d - x^\dag), \la_{n,i}^\d - \la_i \right\r.
\end{align*}
Recall that $a_{n,i}^\d \in \{0, 1\}$. When $a_{n,i}^\d = 0$, we have $\|r_{n,i}^\d\| \le \tau \d_i$ and, by Assumption \ref{ass2} (iv), we obtain
\begin{align*}
\|A_i(x_n^\d - x^\dag)\| \le (1+\eta) \|F_i(x_n^\d) - y_i\| \le (1+\eta) (1+\tau) \d_i.
\end{align*}
When $a_{n, i}^\d = 1$, we have $\|r_{n,i}^\d\| > \tau \d_i$ and hence 
\begin{align*}
\|A_i(x_n^\d - x^\dag)\| \le (1+\eta) \left(\|r_{n,i}^\d\| + \d_i\right) \le (1+\eta) (1+ 1/\tau) \|r_{n,i}^\d\|.  
\end{align*}
Therefore 
\begin{align*}
I_1 \le{}& - \frac{2\gamma_n \Phi_n}{N} \left\l A (x_n^\d - x^\dag), \la_n^\d - \la \right\r
+ \frac{2(1+\eta)(1+\tau)}{N} \gamma_n \Phi_n \sum_{i=1}^N \d_i \|\la_{n,i}^\d - \la_i\| \\
& + \frac{2 (1+\eta)(1+1/\tau) \gamma_n \Phi_n}{N} \sum_{i=1}^N a_{n,i}^\d (1-w_{n,i}^\d) \|r_{n,i}^\d\| \|\la_{n,i}^\d - \la_i\|.
\end{align*}
By using (\ref{smd.111}), the Cauchy-Schwarz inequality, the Young's inequality and $\|\la_{n,i} - \la_i\| \le \|\la_n - \la\|$, we further obtain 
\begin{align}\label{loping_smd.3}
I_1 \le{}& - \frac{2\gamma_n \Phi_n}{N} \left(\Delta_n^\d - \eta_\la\right) 
+ \frac{2(1+\eta)(1+\tau)}{N} \gamma_n \Phi_n \d \|\la_n^\d - \la\| \nonumber \\
& + \frac{2(1+\eta)^2(1+1/\tau)^2 \gamma_n^2 \Phi_n}{N} \sum_{i=1}^N a_{n,i}^\d (1-w_{n,i}^\d) \|r_{n,i}^\d\|^2 \nonumber \\
& + \frac{\Phi_n \|\la_n^\d - \la\|^2}{2N} \sum_{i=1}^N (1-w_{n,i}^\d). 
\end{align}
Next we consider $I_2$. By the Cauchy–Schwarz inequality and the Young's inequality we have 
\begin{align*}
I_2 & \le \frac{(3+\eta)\kappa \gamma_n \Phi_n \|x_n^\d - x^\dag\|}{N} \left(\sum_{i=1}^N a_{n,i}^\d (w_{n,i}^\d)^2 \|r_{n,i}^\d\|^2\right)^{\frac{1}{2}} \left(\sum_{i=1}^N a_{n,i}^\d \|\la_{n,i}^\d - \la_i\|^2\right)^{\frac{1}{2}}  \\
& \le \frac{(3+\eta)\kappa \gamma_n \Phi_n \|x_n^\d - x^\dag\|}{N} \|\la_n^\d - \la\| \left(\sum_{i=1}^N a_{n,i}^\d (w_{n,i}^\d)^2 \|r_{n,i}^\d\|^2\right)^{\frac{1}{2}} \\
& \le \frac{\gamma_n \Phi_n}{N} \left(\sigma \|x_n^\d - x^\dag\|^2 + \frac{(3+\eta)^2 \kappa^2}{4\sigma} \|\la_n^\d - \la\|^2 \sum_{i=1}^N a_{n,i}^\d (w_{n,i}^\d)^2 \|r_{n,i}^\d\|^2 \right)  \\
& \le \frac{\gamma_n \Phi_n}{N} \Delta_n^\d + \frac{(3+\eta)^2 \kappa^2 \gamma_n \Phi_n \|\la_n^\d - \la\|^2}{4\sigma N} \sum_{i=1}^N a_{n,i}^\d (w_{n,i}^\d)^2 \|r_{n,i}^\d\|^2,
\end{align*}
where for the last step we used $\sigma \|x_n^\d - x^\dag\|^2 \le \Delta_n^\d$ which follows from the strong convexity of $\R$. Now we use the elementary inequality 
\begin{align*}
1-e^{-t} = \int_0^t e^{-s} ds \ge t e^{-t} \ge t e^{-2t}, \quad \forall t \ge 0,
\end{align*}
the definition of $w_{n,i}^\d$, and the choice of $\vartheta$ to obtain 
\begin{align*}
\frac{(3+\eta)^2 \kappa^2 \gamma_n }{2\sigma} a_{n,i}^\d \|r_{n,i}^\d\|^2 (w_{n,i}^\d)^2
& = \vartheta c_0 \gamma_n a_{n,i}^\d \|r_{n,i}^\d\|^2 e^{-2\vartheta c_0\gamma_n a_{n,i}^\d \|r_{n,i}^\d\|^2} \\
& \le 1 - e^{-\vartheta c_0\gamma_n a_{n,i}^\d \|r_{n,i}^\d\|^2} \\
& = 1-w_{n,i}^\d.
\end{align*}
Therefore 
\begin{align}\label{loping_smd.4}
I_2 \le \frac{\gamma_n \Phi_n}{N} \Delta_n^\d + \frac{\Phi_n \|\la_n^\d - \la\|^2}{2N} \sum_{i=1}^N (1-w_{n,i}^\d).
\end{align}
Combining (\ref{loping_smd.3}) and (\ref{loping_smd.4}) with (\ref{loping_smd.5}), and using 
\begin{align*}
\frac{(2+(1+\eta)\kappa\rho) \gamma_n \Phi_n}{N} \sum_{i=1}^N w_{n,i}^\d a_{n,i}^\d \d_i \|\la_{n,i}^\d - \la_i\|
\le \frac{(2+(1+\eta)\kappa\rho) \gamma_n \Phi_n}{N} \d \|\la_n^\d - \la\|, 
\end{align*}
we obtain 
\begin{align*}
& \EE_{i_n} \left[\Phi_{n+1} \|\la_{n+1}^\d - \la\|^2\right] \nonumber \\
& \le \Phi_n \|\la_n^\d - \la\| - \frac{\gamma_n \Phi_n}{N} \Delta_n^\d + \frac{2 \gamma_n \Phi_n}{N} \eta_\la \\
& \quad \, \frac{\gamma_n^2 \Phi_n}{N} \sum_{i=1}^N a_{n,i}^\d \left[(1+\kappa\rho)^2 w_{n,i}^\d + 2(1+\eta)^2(1+1/\tau)^2 (1-w_{n,i}^\d)\right] \|r_{n,i}^\d\|^2 \\
& \quad \, + \frac{\gamma_n \Phi_n}{N} \left(2 + (1+\eta)(2+2\tau+\kappa\rho)\right) \d \|\la_n^\d - \la\|.
\end{align*}
Since $0 < w_{n,i}^\d \le 1$, we therefore complete the proof. 
\end{proof}

\begin{lemma}\label{loping_smd.lem4}
Let Assumption \ref{ass2} hold. Consider Algorithm \ref{alg:SMD0} with $\{\gamma_n\}$ satisfying (\ref{sys.2}). Let $\vartheta$ and $\Phi_n$ be as in Lemma \ref{loping_smd.lem3}. Set $\Gamma := e^{\vartheta \sigma \rho^2}$, $\a := C_2 \Gamma \bar\gamma/c_1$ and for any $\la \in Y$ define 
\begin{align*}
\mathcal U_n := \Phi_n \|\la_n^\d-\la\|^2 + \a \Delta_n^\d, \quad n = 0, 1, \cdots. 
\end{align*}
There holds 
\begin{align*}
\EE_{i_n}\left[\mathcal U_{n+1}\right] + \frac{\gamma_n}{N} \Delta_n^\d 
\le \mathcal U_n + \frac{2 \Gamma \gamma_n}{N} \eta_\la + \frac{C_3 \sqrt{\Gamma} \d}{N} \gamma_n \sqrt{\mathcal U_n}.
\end{align*}
\end{lemma}

\begin{proof}
Let $h_n := \gamma_n/N$. It follows from Lemma \ref{loping_smd.lem3} that 
\begin{align*}
\EE_{i_n}\left[\Phi_{n+1}\|\la_{n+1}^\d - \la\|^2 \right]
\le{}& \Phi_n \|\la_n^\d - \la\|^2 - \Phi_n h_n \Delta_n^\d + 2 \Phi_n h_n \eta_\la \\
& + C_2 \Phi_n \gamma_n h_n \sum_{i=1}^N a_{n,i}^\d \|r_{n,i}^\d\|^2 \\
& + C_3 \d h_n \Phi_n \|\la_n^\d - \la\|.
\end{align*}
From Lemma \ref{Land:lem6} we also have 
\begin{align*}
\EE_{i_n}\left[\Delta_{n+1}^\d\right] \le \Delta_n^\d - c_1 h_n \sum_{i=1}^N a_{n,i}^\d \|r_{n,i}^\d\|^2.
\end{align*}
From the above two equations it follows that 
\begin{align*}
\EE_{i_n}\left[\mathcal U_{n+1} \right]
\le{}& \mathcal U_n - \Phi_n h_n \Delta_n^\d + 2 \Phi_n h_n \eta_\la \\
& + h_n \left(C_2 \Phi_n \gamma_n - \a c_1 \right) \sum_{i=1}^N a_{n,i}^\d \|r_{n,i}^\d\|^2 \\
& + C_3 \d h_n \Phi_n \|\la_n^\d - \la\|.
\end{align*}
Note that $C_2 \Phi_n \gamma_n \le C_2 \Gamma \bar\gamma$. Since $\a = C_2 \Gamma \bar\gamma/c_1$, we have $C_2 \Phi_n \gamma_n - \a c_1 \le 0$ and thus  
\begin{align*}
\EE_{i_n}\left[\mathcal U_{n+1} \right] + \Phi_n h_n \Delta_n^\d 
\le \mathcal U_n + 2 \Phi_n h_n \eta_\la + C_3 \d h_n \Phi_n \|\la_n^\d - \la\|.
\end{align*}
Because $1 \le \Phi_n \le \Gamma$, we have 
\begin{align*}
\Phi_n h_n \Delta_n^\d \ge h_n \Delta_n^\d, \qquad 2 \Phi_n h_n \eta_\la \le 2\Gamma h_n \eta_\la.
\end{align*}
Therefore 
\begin{align*}
\EE_{i_n}\left[\mathcal U_{n+1}\right] + h_n \Delta_n^\d 
\le \mathcal U_n + 2 \Gamma h_n \eta_\la + C_3 \d h_n \Phi_n \|\la_n^\d - \la\|.
\end{align*}
Combining this with 
\begin{align*}
C_3 \d h_n \Phi_n \|\la_n^\d - \la\| \le C_3 \d h_n \sqrt{\Gamma \Phi_n \|\la_n^\d-\la\|^2}
\le C_3 \d h_n \sqrt{\Gamma \mathcal U_n},
\end{align*}
we thus complete the proof. 
\end{proof}

According to estimate in Lemma \ref{loping_smd.lem4}, by taking the full expectation and using $\EE[\sqrt{\mathcal U_n}] \le \EE[\mathcal U_n]^{1/2}$ we obtain 
\begin{align*}
\EE\left[\mathcal U_{n+1}\right] + \frac{\gamma_n}{N} \EE\left[\Delta_n^\d\right] 
\le \EE\left[\mathcal U_n\right] + \frac{2 \Gamma \gamma_n}{N} \eta_\la + \frac{C_3 \sqrt{\Gamma} \d \gamma_n}{N} \EE\left[\mathcal U_n\right]^{1/2}.
\end{align*}
Recursively using this equation and noting $\mathcal U_0 \le \Gamma \|\la\|^2 + \a \Delta_0$ gives 
\begin{align*}
\EE\left[\mathcal U_n\right] + \frac{1}{N} \sum_{k=0}^{n-1} \gamma_k \EE\left[\Delta_k^\d\right] 
\le \a \Delta_0 + \Gamma \|\la\|^2 + \frac{2 \Gamma s_n}{N} \eta_\la + \frac{C_3 \sqrt{\Gamma} \d }{N} \sum_{k=0}^{n-1} \gamma_k \EE\left[\mathcal U_k\right]^{1/2}.
\end{align*}
Now we may use the same argument in the proof of Lemma \ref{SMD.lem5} to obtain the following result. 

\begin{proposition}\label{smd.prop1}
Let Assumption \ref{ass2} hold and consider Algorithm \ref{alg:SMD0} with $\{\gamma_n\}$ satisfying (\ref{sys.2}). Then for any $\la \in Y$ there hold
\begin{align*}
\EE\left[\Delta_n^\d\right]
\le C \left( \frac{1}{1+s_n/N} + \frac{\|\la\|^2}{1+s_n/N} + \eta_\la + \frac{s_n}{N} \d^2\right)
\end{align*}
for all $n \ge 0$, where $\eta_\la$ is defined in (\ref{Fenchel-defect}), and $C$ is a constant depending only on $\bar \gamma$, $\sigma$, $L$, $\tau$, $\rho$ and $\Delta_0:= D_\R^{\xi_0}(x^\dag, x_0)$. 
\end{proposition}

The preceding proposition bounds the reconstruction error in terms of an arbitrary dual comparison element and its associated Fenchel defect. Based on this result, we finally give convergence rates for Algorithm \ref{alg:SMD0} under local variational source conditions for nonlinear problems. 

\begin{theorem}\label{thm:LSMD}
Let Assumption \ref{ass2} hold. Consider Algorithm \ref{alg:SMD0} and assume that (\ref{sys.2}) and (\ref{SMD.44}) are fulfilled. Let $n_\d$ be the random integer determined by (\ref{LSMD-stop}), which is finite almost surely. If there exist $M \ge 0$ and $0< q\le 1$ such that 
\begin{align}\label{lvsc}
\R(x^\dag) - \R(x) - \l \xi_0, x^\dag - x\r \le M \|F(x)-y\|^q
\end{align}
holds for every $x\in B_{\rho}(x_0)\cap {\rm dom}(\R)$, then  
\begin{align*}
\EE\left[D_\R^{\xi_{n_\d}^\d}(x^\dag, x_{n_\d}^\d) \right] \le C \d^q \quad \mbox{and} \quad 
\EE\left[\|x_{n_\d}^\d - x^\dag\|^2 \right] \le C \d^{q},
\end{align*}
where $C$ is a constant depending only on $\bar \gamma$, $\sigma$, $L$, $\tau$, $\rho$ and $\Delta_0:= D_\R^{\xi_0}(x^\dag, x_0)$.
\end{theorem}

\begin{proof}
By Assumption \ref{ass2} (iv) we have $\|F(x) - y\|\le \frac{1}{1-\eta} \|A(x - x^\dag)\|$ for $x \in B_{2\rho}(x_0)$. Thus, it follows from the variational source condition (\ref{lvsc}) that 
\begin{align}\label{smd.101}
D_\R^{\xi_0}(x^\dag, x_0) - D_\R^{\xi_0}(x, x_0) \le M \|F(x) - y\|^q \le \frac{M}{(1-\eta)^q} \|A(x -x^\dag)\|^q
\end{align}
for all $x \in B_\rho(x_0) \cap {\rm dom}(\R)$. By virtue of (ii) and (iii) in Assumption \ref{ass2}, we have for $x \not\in B_\rho(x_0)$ that 
\begin{align*}
D_\R^{\xi_0}(x, x_0) \ge \sigma \|x-x_0\|^2\ge \sigma \rho^2 \ge D_\R^{\xi_0}(x^\dag, x_0).
\end{align*}
Therefore (\ref{smd.101}) holds for all $x \in {\rm dom}(\R)$. Consequently 
\begin{align*}
\R(x^\dag) - \R(x) - \l \xi_0, x^\dag - x\r \le \frac{M}{(1-\eta)^q} \|A(x-x^\dag)\|^q, \quad \forall x \in {\rm dom}(\R).
\end{align*}
Based on this inequality, Lemma \ref{Land:lem6} and Proposition \ref{smd.prop1}, we may use the similar arguments in the proofs of Lemma \ref{SMD.lem6} and Theorem \ref{SMD.thm5} to conclude the desired result. 
\end{proof}

\begin{remark}
The stochastic gradient method for nonlinear ill-posed inverse problems in Hilbert spaces was studied in \cite{JZZ2020}. Under the diminishing step size $\gamma_n = \gamma_0 (n+1)^{-\a}$ with $\a\in (0,1)$ being sufficiently close to $1$, \cite[Theorem 2.5]{JZZ2020} shows that, under the H\"{o}lder source condition 
$$
x_0 -x^\dag = (A^* A)^\nu \omega
$$ 
for some $0< \nu \le 1/2$ and $\omega \in X$, the convergence rate
$$
\EE\left[\|x_{\hat n_\d}^\d - x^\dag\|^2\right] = O\left(\d^{\frac{4\nu}{1+2\nu}}\right). 
$$
holds for a suitably chosen \emph{a priori} deterministic stopping index $\hat n_\d$. In addition to the standard assumptions comparable to Assumption \ref{ass2}, the analysis in \cite{JZZ2020} relies on the technical conditions
\begin{align*}
\EE\left[\|(I - R_{z_t}) G(x_n^\d)\|^2\right]^{1/2} 
&\le C \EE\left[\|x_n^\d-x^\dag\|^2\right]^{\theta/2}\EE\left[\|G(x_n^\d)\|^2\right]^{1/2}, \\
\EE\left[\|(I - R_{z_t}^*) G(x_n^\d)\|^2\right]^{1/2} 
&\le C \EE\left[\|x_n^\d-x^\dag\|^2\right]^{\theta/2}\EE\left[\|G(x_n^\d)\|^2\right]^{1/2}
\end{align*}
for some $\theta\in (0, 1]$; see \cite[Assumption 2.4]{JZZ2020} for the definitions of $R_{z_t}$ and $G(x)$. Unfortunately, no verification of these technical conditions is currently available. 

In contrast, our analysis does not require such additional technical assumptions and yields order-optimal convergence rates under the stated structural assumptions. Moreover, our results go beyond the a priori framework by establishing order-optimal rates for the iterate selected by an \emph{a posteriori} stopping rule.
\end{remark}

\section{\bf Numerical results}

In this section, we present numerical results to evaluate the performance of Algorithms \ref{alg:SMD2} and \ref{alg:SMD0}. Both algorithms sample the index $i_n$ uniformly from $\{1,\cdots,N\}$. If $i_n\not\in I_n^\d$, the corresponding equation has already been verified at the current iterate, since the verification set is reset after every update. Thus, $a_{n,i_n}=0$, and the selection leaves both the iterates and the verification set unchanged:
\begin{align*}
\xi_{n+1}^\d = \xi_n^\d,\qquad x_{n+1}^\d = x_n^\d,\qquad I_{n+1}^\d = I_n^\d.
\end{align*}
To avoid these redundant iterations, we modify the implementation by sampling $i_n$ uniformly from the current verification set $I_n^\d$. We use this modification in all numerical experiments below. Given the extensive numerical studies of stochastic methods already available (\cite{GJ2026,HJLZ2025,JL2019,JLT2026,JLZ2023}), our experiments focus on assessing the convergence rates of the reconstructions produced at the verification-based stopping indices.

\begin{example}
We consider the Fredholm integral equation of the first kind
\begin{align*}
(Ax)(s) := \int_0^1 k(s,t) x(t) dt = y(s), \qquad 0\le s\le 1,
\end{align*}
with the Gaussian kernel $k(s,t) := e^{-\frac{(s-t)^2}{2a^2}}$, where $a = 0.05$. It is clear that $A: L^2(0,1) \to L^2(0,1)$ is a compact linear operator. We assume the sought solution is sparse and use the sparsity-promoting functional
\begin{align*}
\R(x) = \frac{1}{2} \|x\|_{L^2(0,1)}^2 + \beta \|x\|_{L^1(0,1)},
\end{align*}
where $\beta = 1$. Then, for any given $\xi \in L^2(0,1)$ we have 
\begin{align*}
x_\xi := \arg\min_{x\in L^2(0,1)} \left\{\R(x) - \l \xi, x\r\right\} = S_\beta(\xi),
\end{align*}
where $S_\beta$ is the soft-thresholding operator
\begin{align*}
S_\beta(z) = {\rm sign}(z) \max\{|z|-\beta, 0\}.
\end{align*}
To verify the theoretical convergence rate result for Algorithm \ref{alg:SMD2}, we generate the sparse exact solution $x^\dag$ by $x^\dag = S_\beta(\xi_0 + A^*\la^\dag)$, where 
\begin{align*}
\la^\dag(s) ={}& 21e^{-((s-0.18)/0.035)^2} - 30 e^{-((s-0.42)/0.05)^2} + 24e^{-((s-0.68)/0.04)^2} \\
& - 18e^{-((s-0.86)/0.03)^2},
\end{align*}
which satisfies the variational source condition with $q=1$. 

\begin{table}[ht]
\centering
\caption{Numerical results for the first kinf Fredholm integral equation with a sparse solution.}
\label{tab:sparse-fredholm}
\begin{tabular}{ccccc}
\hline
\(\d_{\rm rel}\) & $\d$ & $\EE\left[\|x_{n_\d}^\d-x^\dag\|^2\right]$ & ratio & mean $n_\d$ \\
\hline
\(5.0\times10^{-1}\) & \(1.5487\times10^{-2}\) & \(6.2952\times10^{-3}\) & 0.4065 & \(87.1\) \\
\(5.0\times10^{-2}\) & \(1.5487\times10^{-3}\) & \(7.4495\times10^{-4}\) & 0.4810 & \(167.3\) \\
\(5.0\times10^{-3}\) & \(1.5487\times10^{-4}\) & \(7.3912\times10^{-5}\) & 0.4773 & \(1620.8\) \\
\(5.0\times10^{-4}\) & \(1.5487\times10^{-5}\) & \(4.7511\times10^{-6}\) & 0.3068 & \(10526.6\) \\
\(5.0\times10^{-5}\) & \(1.5487\times10^{-6}\) & \(1.2378\times10^{-7}\) & 0.0799 & \(30813.0\) \\
\(5.0\times10^{-6}\) & \(1.5487\times10^{-7}\) & \(1.3862\times10^{-8}\) & 0.0895 & \(113364.8\) \\
\hline
\end{tabular}
\end{table}

For numerical implementation, we discretize the equation by using $n_x = 600$ midpoint quadrature points for the unknown and $n_y = 800$ midpoint quadrature points for the data. The data are divided into $N = 40$ interlaced blocks. For a prescribed relative noise level \(\delta_{\rm rel}\), the total noise level is chosen as $\d =
\d_{\rm rel} \|y\|$, and balanced component noise levels are defined by $\d_i = \d/\sqrt{N}$. For each block, an independent Gaussian noise vector is generated and rescaled so that
\begin{align*}
\|y_i^\d - y_i\| = \d_i.
\end{align*}
We apply Algorithm \ref{alg:SMD2} with the modification that, at iteration $n$, the index $i_n$ is sampled uniformly from the current verification set $I_n^\d$. The iteration terminates when the verification set becomes empty. In the experiment we use the initial guess $\xi_0 = 0$ and $x_0 = 0$, the parameter $\tau=1.2$ and the constant step size
\begin{align*}
\gamma_n \equiv \gamma = \frac{1.98}{L} \left(1-\frac{1}{\tau}\right), \qquad
L = \max_{1\le i\le N} \|A_i\|^2.
\end{align*}
The theoretical convergence result in Theorem \ref{SMD.thm5} for $q=1$ predicts
\begin{align*}
\EE\left[\|x_{n_\delta}^\delta-x^\dagger\|_{L^2}^2\right] = O(\d).
\end{align*}
To examine this rate numerically, we consider 6 relative noise levels ranging from $5 \times 10^{-1}$ to $5 \times10^{-6}$. The expectation is approximated by Monte Carlo averaging over 20 independent runs. The results are reported in Table \ref{tab:sparse-fredholm}, including the mean squared reconstruction error, the normalized quantity 
\begin{align*}
\texttt{ratio} := \EE\left[\|x_{n_\d}^\d - x^\dag\|_{L^2}^2\right]/\d,
\end{align*}
and the mean stopping index. From Table \ref{tab:sparse-fredholm} we can see that the quantity ``\texttt{ratio}" remains bounded throughout the entire noise range. This is consistent with the theoretically predicted rate estimate. The stopping index increases substantially as the noise decreases, reflecting the increasingly stringent discrepancy thresholds in the verification procedure.
\end{example}

\begin{example}
We next consider the identification of a comment reaction coefficient $c$ from measurements of the solutions to a family of boundary value problems. In the $i$-th experiment, we prescribe a source $f_i \in L^2(0,1)$ and boundary values $g_{i,0}, g_{i,1} \in \mathbb R$. For a coefficient $c$, the corresponding state $u_i(c)$ is determined by
\begin{align*}
\begin{cases}
-u_i''(t) + c(t) u_i(t) = f_i(t), & 0<t<1,\\
u_i(0) = g_{i,0}, \qquad u_i(1) = g_{i,1},
\end{cases}
\qquad i = 1, \cdots, N.
\end{align*}
We take the coefficient and observation spaces to be $X=L^2(0,1)$ and $Y_i=L^2(0,1)$. According to the theory of elliptic equation, there exists a constant $\ep_0 > 0$ such that for every  
\begin{align*}
c \in \mathcal D := \{c\in L^2(0,1): \|c - \hat c\|_{L^2} < \ep_0 \text{ for some } \hat c \ge 0\text{ a.e.}\},
\end{align*}
the boundary value problem has a unique solution \(u_i(c)\in H^2(0,1)\). We therefore define the \(i\)-th forward operator by
\begin{align*}
F_i: \mathcal D \subset X \to Y_i, \qquad F_i(c) = u_i(c).
\end{align*}
Let $c^\dag\in\mathcal D$ denote the coefficient to be recovered. In each experiment, we observe the corresponding state on \((0,1)\). The exact observations are $y_i = F_i(c^\dag)$, $i = 1, \cdots, N$, and the available noisy measurements $y_i^\d \in Y_i$ satisfy
\begin{align*}
\|y_i^\d - y_i\|_{L^2(0,1)} \le \d_i.
\end{align*}
The inverse problem is to recover the common coefficient $c^\dag$ from these measurements, or equivalently to solve the nonlinear system
\begin{align*}
F_i(c) = y_i, \qquad i = 1, \cdots, N,
\end{align*}
using the noisy data $y_i^\d$.

The use of several experiments allows us to collect solution profiles generated by different sources or boundary excitations. These profiles can provide complementary information about the coefficient, particularly in regions where the state from one experiment is small and hence responds weakly to changes in \(c\). Here \(N\) is simply the number of available experiments.

It is known that, for each $c \in \mathcal D$, every forward operator $F_i$ is Fr\'{e}chet differentiable at $c$ with Fr\'{e}chet derivative and the adjoint given by 
\begin{align*}
F_i'(c) h = v_i \quad \mbox{and} \quad F_i'(c)^* r_i = - u_i(c) p_i
\end{align*}
for $h \in X$ and $r_i \in Y_i$, where $v_i$ and $p_i$ are characterized by
\[
\begin{cases}
-v_i''+cv_i=-h\,u_i(c) \quad \mbox{in } (0, 1),\\
v_i(0)=v_i(1)=0,
\end{cases}
\qquad 
\begin{cases}
-p_i''+cp_i=r_i \quad \mbox{in } (0,1),\\
p_i(0)=p_i(1)=0.
\end{cases}
\]

We assume that the sought coefficient $c^\dag$ is nonnegative. To incorporate the nonnegativity constraint into the reconstruction, we use
\begin{align*}
\R(c) = \frac{1}{2} \|c\|_{L^2(0,1)}^2 + \iota_{\mathcal C}(c),
\end{align*}
where $\iota_{\mathcal C}$ is the indicator function of the set $\mathcal C := \{c \in X: c \ge 0 \text{ a.e.}\}$.  It is easy to see that, for any $\xi \in L^2(0,1)$ there holds 
\begin{align*}
\arg\min_{x\in L^2(0,1)} \left\{\R(x) - \l \xi, x\r\right\} = \max\{\xi, 0\}.
\end{align*}

In our numerical experiment, we choose the exact coefficient
\begin{align*}
c^\dag(t) = \max\{0,t^2(8-10t)^3\}
\end{align*}
which is nonnegative, compactly supported in \([0,0.8]\), and belongs to \(H^2(0,1)\). The exact states are prescribed by
\begin{align*}
u_i^\dag(t) = 1 + \ep_i\cos(i\pi t), \qquad i = 1, \cdots, N,
\end{align*}
with $\ep_i = 0.05 + 0.15\frac{i-1}{N-1}$. Hence $u_i^\dag(t) \ge 0.8$ and the boundary data are 
\begin{align*}
g_{i,0} = 1 + \ep_i, \qquad g_{i,1} = 1 + (-1)^i \ep_i.
\end{align*}
The forcing terms are chosen as $f_i = -(u_i^\dag)'' + c^\dag u_i^\dag$ so that $F_i(c^\dag) = u_i^\dag$. According to \cite[Example 4.2]{HNS95}, every $F_i$ satisfies Assumption \ref{ass2} around $c^\dag$. 

The chosen exact solution satisfies the variational source condition with $q=1$. Indeed, let
\begin{align*}
p_i^\dag = - \frac{c^\dag}{N u_i^\dag}, \qquad i = 1, \cdots, N.
\end{align*}
Then $p_i^\dag \in H^2(0,1)\cap H_0^1(0,1)$. Set $\la_i^\dag = -(p_i^\dag)'' + c^\dag p_i^\dag$. Then $\la_i^\dag \in L^2(0,1)$ and  
\begin{align*}
F_i'(c^\dag)^*\la_i^\dag = -u_i^\dag p_i^\dag = \frac{1}{N} c^\dag.
\end{align*}
Summing over $i$ gives $F'(c^\dag)^*\la^\dag = c^\dag$. With $c_0 = \xi_0 = 0$, this yields $F'(c^\dag)^*\la^\dag
\in \p\R(c^\dag)$, and hence the variational source condition holds with $q=1$.

We generate the noisy data as follows. For a prescribed relative noise level $\d_{\rm rel}$ we set $\d = \d_{\rm rel} \|y\|$, where $\|y\|^2 = \sum_{i=1}^N \|y_i\|_{L^2}^2$. For each component, an independent Gaussian noise vector is generated and rescaled so that 
\begin{align*}
\|y_i^\delta-y_i\|_{L^2} = \d/\sqrt{N}.
\end{align*}
We apply Algorithm \ref{alg:SMD0} with the modification that, at each iteration, the index $i_n$ is sampled uniformly from the current verification set $I_n^\d$. The method terminates when the verification set becomes empty. The theoretical result in Theorem \ref{thm:LSMD} for $q=1$ predicts
\begin{align*}
\EE\left[\|c_{n_\d}^\d - c^\dag\|_{L^2}^2\right] = O(\d).
\end{align*}

\begin{table}[ht]
\centering
\caption{Numerical results for the nonlinear reaction-coefficient identification problem.}
\label{tab:nonlinear}
\begin{tabular}{ccccc}
\hline
\(\d_{\rm rel}\) & \(\d\) & \(\EE\left[\|c_{n_\d}^\d-c^\dag\|_{L^2}^2\right]\) & ratio & mean \(n_\d\) \\
\hline
\(5.0\times10^{-2}\) & \(3.1717\times10^{-1}\) & \(2.0179\)              & 6.3622 & 131.8    \\
\(1.0\times10^{-2}\) & \(6.3433\times10^{-2}\) & \(3.3000\times10^{-1}\) & 5.2022 & 276.0 \\
\(5.0\times10^{-3}\) & \(3.1717\times10^{-2}\) & \(2.3248\times10^{-1}\) & 7.3300 & 434.2  \\
\(1.0\times10^{-3}\) & \(6.3433\times10^{-3}\) & \(4.0934\times10^{-2}\) & 6.4531 & 6059.2 \\
\(5.0\times10^{-4}\) & \(3.1717\times10^{-3}\) & \(1.9438\times10^{-2}\) & 6.1287 & 12480.0 \\
\(1.0\times10^{-4}\) & \(6.3433\times10^{-4}\) & \(2.8991\times10^{-3}\) & 4.5703 & 41563.1 \\
\(5.0\times10^{-5}\) & \(3.1717\times10^{-4}\) & \(1.2893\times10^{-3}\) & 4.0652 & 67420.7 \\
\(1.0\times10^{-5}\) & \(6.3433\times10^{-5}\) & \(3.0006\times10^{-4}\) & 4.7303 & 194296.8 \\ \hline
\end{tabular}
\end{table}

In the numerical experiment we use $N= 40$ and $\tau = 1.2$. We discretize the interval $(0,1)$ by $m=400$ interior grid points and approximate $-d^2/dt^2$ by the standard second-order centered difference scheme. We use the constant step size $\gamma_n \equiv \gamma = 40$. The numerical results are summarized in Table~\ref{tab:nonlinear}, where the mean squared reconstruction error $\EE\left[\|c_{n_\d}^\d-c^\dag\|_{L^2}^2\right]$ is calculated approximately by Monte Carlo averaging over 20 independent runs, and ``\texttt{ratio}" denotes $\EE \left[\|c_{n_\d}^\d - c^\dag\|_{L^2}^2 \right]/\d$.  We consider 8 relative noise levels ranging from $5\times 10^{-2}$ to $1 \times 10^{-5}$ and the values of \texttt{ratio} are between approximately 4.06 and 7.33. This behavior is consistent with the theoretically predicted estimate for $q = 1$. 
\end{example}

\section{\bf Conclusion}

We have developed a quantitative regularization theory for the loping stochastic mirror descent method applied to systems of linear and nonlinear ill-posed equations with noisy data. The method combines stochastic component updates with a verification-set mechanism that provides an implementable \emph{a posteriori} stopping rule without requiring scheduled evaluations of the full residual.

A central point of the analysis is that the convergence rates are established for the stopped iterate actually returned by the algorithm. Although the stopping index is random and is strongly coupled to the stochastic trajectory, its distribution does not have to be estimated. The pathwise monotonicity of the Bregman error, together with stationarity of the iteration after successful verification, allows the stopped error to be controlled by deterministic-horizon estimates.

For linear problems, we derived deterministic-index estimates involving arbitrary dual comparison elements and their Fenchel defects. This formulation makes it possible to exploit variational source conditions without requiring an exact source representation. Optimizing the resulting approximation--stability balance yields the order-optimal convergence rates. 

For nonlinear systems, the principal difficulty is the appearance of perturbation terms coupling the primal error, component residuals, and auxiliary dual variables. We introduced an exponential-weighted energy argument that converts the Bregman descent into the damping needed to control these nonlinear interactions. Under standard nonlinearity conditions on forward operators and local variational source conditions on sought solutions, the same order-optimal convergence rates are obtained. 

Several questions remain open. First, it would be interesting to strengthen the present rate results in expectation to obtain almost sure convergence rates for the discrepancy-stopped iterates. Such results would require a finer pathwise analysis of the stochastic errors and of the random stopping index, possibly using martingale convergence and concentration techniques. Second, adaptive step-size strategies have already been developed, and their convergence has been established (\cite{GJ2026,JLZ2023,HJLZ2025}). A natural next question is whether the convergence-rate guarantees obtained here extend to such strategies. The main difficulty is that adaptive step sizes depend on the random iteration history and may also depend on the currently sampled equation, coupling them with the residuals and reconstruction errors. The resulting accumulated step size is therefore random, which complicates the conditional-expectation estimates and their transfer to the stopping index. Establishing order-optimal rates for the actual iterates selected by the discrepancy principle under adaptive step-size rules requires quantitative control of these dependencies.

\bibliographystyle{amsplain}

\end{document}